\documentclass[a4paper,12pt]{article}
\usepackage{cmap}
\usepackage[T1]{fontenc}
\usepackage[utf8]{inputenc}
\usepackage[a4paper]{geometry}
\usepackage{enumitem}
\usepackage[backend=bibtex,maxnames=10,defernumbers=true,giveninits=true,sortcites=true,isbn=false]{biblatex}
\usepackage[unicode,hyperfootnotes=false,pdftex]{hyperref}
\usepackage{bookmark}
\usepackage{amsthm,amssymb}
\usepackage{mathtools}
\usepackage{microtype}
\usepackage{hyperxmp}
\usepackage{csquotes}
\usepackage{float}
\usepackage[english]{babel}
\usepackage{tikz}
\usepackage{stmaryrd}
\usepackage{mathabx}
\usepackage{tikz-cd}

\DeclareFontFamily{U}{BOONDOX-calo}{\skewchar\font=45 }
\DeclareFontShape{U}{BOONDOX-calo}{m}{n}{
  <-> s*[1.05] BOONDOX-r-calo}{}
\DeclareFontShape{U}{BOONDOX-calo}{b}{n}{
  <-> s*[1.05] BOONDOX-b-calo}{}
\DeclareMathAlphabet{\mathcalboondox}{U}{BOONDOX-calo}{m}{n}
\SetMathAlphabet{\mathcalboondox}{bold}{U}{BOONDOX-calo}{b}{n}
\DeclareMathAlphabet{\mathbcalboondox}{U}{BOONDOX-calo}{b}{n}

\hypersetup{%
    pdftitle={Beyond the delta method},
    pdfauthor={Antoine Lejay; Sara Mazzonetto},
    pdfdate={2024-01-31},
    pdfkeywords={parametric estimation; Taylor expansion; hypothesis testing; delta method; maximum likelihood estimation; moments method; Implicit Function Theorem},
    pdflang={en},
    pdfmetalang={en},
    pdfsubject = {    We give an asymptotic development of the maximum likelihood estimator~(MLE), 
    or any other estimator defined implicitly, in a way which involves
    the limiting behavior of the score and its higher-order derivatives. 
    This development, which is explicitly computable, 
    gives some insights about the non-asymptotic behavior of the renormalized MLE and 
    its departure from its limit. 
    We highlight that the results hold whenever the score and its derivative converge, including to non Gaussian limits.
    Our approach is based on an asymptotic implicit function theorem, inspired from perturbative approaches.}
}

\SetEnumitemKey{thm}{nolistsep,topsep=-\parskip,label={\textnormal{(\roman*)}}}
\SetEnumitemKey{compact}{noitemsep,topsep=-\parskip,leftmargin=0em,labelindent=2em}

\newcommand{\upflat}{\tikz \draw[thick] (0,0) -- (3pt,3pt) -- (9pt,3pt);}
\newcommand{\zig}{\tikz \draw[thick] (0,0) -- (3pt,3pt) -- (6pt,0pt) -- (9pt,3pt);}

\newcommand{\dummysymbol}[1][3mm]{%
    \begin{tikzpicture}[scale=0.3]
	\fill (0,0) -- (0,-#1) arc (-90:0:#1) -- (0,0);
	\fill (0,0) -- (0,#1) arc (90:180:#1) -- (0,0);
	\draw (0,0) circle[radius=#1];
    \end{tikzpicture}
}

\newcommand{\dummy}{\dummysymbol}
\newcommand{\dummysmall}{\dummysymbol[2.3mm]} % small size for script size

\newcommand{\cC}{\mathcalboondox{C}}
\newcommand{\ce}{\mathcalboondox{e}}
\newcommand{\cN}{\mathcalboondox{N}}
\newcommand{\cF}{\mathcalboondox{F}}
\newcommand{\bS}{\mathbf{S}}
\newcommand{\cI}{\mathcalboondox{I}}
\newcommand{\cT}{\mathcalboondox{T}}
\newcommand{\cA}{\mathcalboondox{A}}
\newcommand{\cE}{\mathcalboondox{E}}
\newcommand{\cL}{\mathcalboondox{L}}
\newcommand{\sX}{\mathsf{X}}
\newcommand{\cH}{\mathcalboondox{H}}
\newcommand{\fp}{\mathfrak{p}}
\newcommand{\fs}{\mathfrak{s}}

\newcommand{\cvproba}[1][]{\xrightarrow[#1]{\text{prob.}}}
\newcommand{\cv}[1][]{\xrightarrow[#1]{}}
\newcommand{\cvprobat}[1][]{\xrightarrow[#1]{\mathbb{P}_\theta}}
\newcommand{\cvdist}[1][]{\xrightarrow[#1]{\text{dist.}}}

\newcommand{\RR}{\mathbb{R}}
\newcommand{\II}{\mathbb{I}}
\newcommand{\ZZ}{\mathbb{Z}}
\newcommand{\EE}{\mathbb{E}}
\newcommand{\PP}{\mathbb{P}}
\newcommand{\KK}{\mathbb{K}}
\newcommand{\NN}{\mathbb{N}}

\newcommand{\bbeta}{\boldsymbol{\beta}}
\newcommand{\bmu}{\boldsymbol{\mu}}
\newcommand{\bdelta}{\boldsymbol{\delta}}
\newcommand{\bgamma}{\boldsymbol{\gamma}}
\newcommand{\balpha}{\boldsymbol{\alpha}}

\newcommand{\dd}{\mathrm{d}}
\newcommand{\vd}{\,\mathrm{d}}
\newcommand{\ba}{\mathbf{a}}

\newcommand{\mlen}{\theta(n)}
\DeclareMathOperator{\Var}{Var}
\newcommand{\dD}{\mathrm{D}}
\newcommand{\pP}[1]{P_{#1}}
\DeclareMathOperator{\grandO}{O}

\newcommand\given{\nonscript\:\delimsize\vert\nonscript\:\mathopen{}} 
\newcommand\SetSymbol[1][]{\nonscript\:#1\vert\nonscript\:\mathopen{}\allowbreak}
\DeclarePairedDelimiterX\Set[1]\{\}{%
  \renewcommand\given{\SetSymbol[\delimsize]}#1}
\DeclarePairedDelimiterX\Prb[1]{[}{]}{%
  \renewcommand\given{\SetSymbol[\delimsize]}#1}
\DeclarePairedDelimiterX\Paren[1](){#1}
\DeclarePairedDelimiterX\dParen[1]{\delimsize(}{\delimsize)}{#1}
\DeclarePairedDelimiterX\dbra[1]{[\delimsize[}{]\delimsize]}{#1}
\DeclarePairedDelimiter{\abs}{|}{|}
\DeclarePairedDelimiter{\bra}{[}{]}
\DeclarePairedDelimiterX\coef[1]\{\}{#1}

\newtheorem{theorem}{Theorem}
\newtheorem{proposition}{Proposition}
\newtheorem{lemma}{Lemma}
\newtheorem{corollary}{Corollary}

\theoremstyle{definition}
\newtheorem{definition}{Definition}
\newtheorem{example}{Example}
\newtheorem{notation}{Notation}

\theoremstyle{remark}
\newtheorem{hypothesis}{Hypothesis}
\newtheorem{remark}{Remark}

\begin{document}

\title{Beyond the delta method}

\author{Antoine Lejay\footnote{Universit\'e de Lorraine, CNRS, IECL, Inria, F-54000 Nancy, France;
  ORCID: \href{https://orcid.org/0000-0003-0406-9550}{0000-0003-0406-9550}; \texttt{Antoine.Lejay@univ-lorraine.fr}}
    \and Sara Mazzonetto\footnote{Universit\'e de Lorraine, CNRS, IECL, Inria, F-54000 Nancy, France; 
	ORCID: \href{https://orcid.org/0000-0001-6187-2716}{0000-0001-6187-2716};
    \texttt{Sara.Mazzonetto@univ-lorraine.fr}}
}

\date{\today}

\maketitle

\begin{abstract}
    We give an asymptotic development of the maximum likelihood estimator~(MLE), 
    or any other estimator defined implicitly, in a way which involves
    the limiting behavior of the score and its higher-order derivatives. 
    This development, which is explicitly computable, 
    gives some insights about the non-asymptotic behavior of the renormalized MLE and 
    its departure from its limit. 
    We highlight that the results hold whenever the score and its derivative converge, including to non Gaussian limits.
    Our approach is based on an asymptotic implicit function theorem, inspired from perturbative approaches.
\end{abstract}

\begin{quote}
    \small
\textbf{MSC(2020) Classification:} Primary 62F12; Secondary 62F03, 62M02, 26B10.

\textbf{Keywords:} parametric estimation; Taylor expansion; hypothesis testing; delta method; maximum likelihood estimation; moments method; Implicit Function Theorem.
\end{quote}

\section{Introduction}

The traditional approach of dealing with parametric estimation consists 
in defining an estimator as the root or minimizer of a function 
and then to study its asymptotic properties as the size of the observations
increases. This is the case for the Maximum Likelihood Estimation (MLE, \cite{pawitan2001all})
and its variants (pseudo- or quasi- likelihood estimation), 
or the Generalized Method of Moments which are used in a broad 
variety of problems \cite{hall2005generalized, Newey_2004}. A typical statement is that
the estimator is consistent and asymptotically (mixed) normal: The 
properly renormalized error between the estimator and the unknown value
converges in distribution to a (mixed) Gaussian random variable. 

The limiting distribution may be used to quantify the quality of the 
estimation, and to construct for example confidence intervals. As 
pointed by D.A. Sprott \cite{sprott75a,sprott73a}, the quality 
of the approximation depends more on the distance between the law of the estimator and the limiting distribution
than on the size itself. Various tools have been introduced to improve 
the estimator for small or average sample size 
or to quantify its departure from its limiting distribution.
Among them, let us cite 
Edgeworth expansions~\cite{skovgaard},
Berry-Esséen inequality~\cite{bishwal11,wang,mishra},
Stein method~\cite{Anastasiou2017,Anastasiou2017b}, 
jacknife~\cite{Does1988},
bootstrap~\cite{EFRON1985}, 
\textquote{magic formula}~\cite{MR712023},
Bartlett corrections \cite{MR56889,MR3289997}, 
and so on. Some of these corrections use considerations on the
density of the estimator, others are pointwise ones.

We look for an expansion of the univariate root~$\theta(n)$ of $F_n(\theta(n))=0$
---~think of the MLE~--- 
of the form
\begin{equation*}
    \theta(n)=\theta+H(n,G_n)
\end{equation*}
where
$\theta$ is the true, unknown parameter, 
$G_n$ is asymptotical pivotal converging to some $G$ (think of the standard Gaussian random variable
coming from a Central Limit Theorem) and $x\mapsto H(n,x)$ is a non-linear, 
random function of the kind
\begin{equation}
    \label{eq:intro:8}
    H(n,x)=\sum_{k=1}^{p} \alpha_k(n)\varphi(n)^k x^k+\grandO(\varphi(n)^{p+1}).
\end{equation}
In the latter expansion, $\varphi(n)$ decreases to $0$ as $n\to\infty$ and 
\begin{equation}
    \label{eq:intro:9}
    \Set{\alpha_k(n)}_{k=1}^p \xrightarrow[n\to\infty]{\text{law or probability}}
    \Set{\alpha_k}_{k=1}^p
\end{equation}
for some $\alpha_k$'s.
Of course, the $\alpha_k$'s and the 
$\alpha_k(n)$'s may depend on $\theta$. The expansions may have different
forms for some particular values of $\theta$ for which some of the $\alpha_k$
vanishes (which we call a \emph{phase transition}).
We highlight that our approach is pointwise, and the order $p$ is related to the regularity of $F_n$.
Note that the expansion of $\theta(n)$ is such that 
\begin{equation*}
    \varphi(n)^{-1}(\theta(n)-\theta)\cvdist[n\to\infty] c_1G.
\end{equation*}
The main result is provided in Theorem~\ref{th:random} in Section~\ref{sec:IFT}: 
An expansion of $\theta(n)$ as just described.

The family of coefficients $\Set{\alpha_k(n)}_{k=1}^p$ in \eqref{eq:intro:8} is not unique. 
However, in Section~\ref{sec:alpha:construction}, given a function~$F_n$ and knowing the asymptotic behavior 
of $F_n$ and its high order derivatives $\dD^k F_n$, we give a procedure to compute
automatically the $\alpha_k(n)$'s satisfying~\eqref{eq:intro:9}
from $F_n$ and its derivatives. 

The expansion we provide may remind the reader of the Edgeworth expansions in~\cite{skovgaard} for the law of the MLE, but our result is pathwise.

Our approach is based on an 
asymptotic implicit function theorem (see Theorem~\ref{thm:1} below), inspired from perturbative approaches. 
Performing asymptotic expansions has a long history with a broad range of
applications, notably to study perturbations of (partial) differential
equations and dynamical
systems~\cite{zbMATH01956296,zbMATH00041890,zbMATH03439409} or in change of
scale analysis \cite{zbMATH05969631}. 
Although in one case, we recover the Lagrange inversion formula 
for power series \cite{gessel,sokal,Furstenberg1967}, our starting point is different
and seems not to have been considered in a general form but mostly for 
particular cases. 
Moreover, we also consider functions that are not necessarily analytic.
Our approach seems to be new because it assumes 
weak general assumptions on the regularity of $F_n$ 
and takes into account different possible asymptotic behaviors of $\dD^k F_n$
(for statistical applications, think that $\dD^k F_n$ may converge
thanks to the Law of Large Numbers or to the Central Limit Theorem, hence with 
different rates). 

Theorem~\ref{th:random} is a sort of generalization of the Delta method \cite{MR1652247} in the sense
that knowing the asymptotic behavior of the score and its derivative,
we get the asymptotic behavior of a function of it, in particular the estimator.
Besides, we do not restrict ourselves to a first order linearization.

In Section~\ref{sec:change-var} we study the effects on our procedure 
of change of variable for reparametrization
of the space of parameters. Related to the problem of reparametrization, 
changing $F_n$ to $F_n/\Phi$ for a suitable $\Phi$
does not change the root~$\theta(n)$ but leads to a different expansion 
in~\eqref{eq:intro:8} which is sometimes more convenient
(see Section~\ref{sec:change-scale} and Section~\ref{sec:exp-fam:natural}).

We illustrate our results through several known examples:
exponential family (Section~\ref{sec:exp-fam})
and parametric estimation for stochastic 
processes (Section~\ref{sec:diffusion}).
A far less trivial application is provided in~\cite{lm22} 
where we apply our procedure to the Skew Brownian 
motion obtaining new expansions of the MLE estimator of the skewness parameter
and showing that it is asymptotically mixed normal. 
In this context, when Skew Brownian motion is a Brownian motion, we recover the expansion given in \cite{lejay2014} 
and give an alternative one obtained via a finer 
expansion for some specific asymptotics of the $\dD^k F_n$.
We refer to this situation of getting two expansions as \emph{phase transition}, 
which we explain in this paper on the binomial model (see Section~\ref{sec:symmetric}).

The examples we mention show that we are not restricted to independent samples
nor to asymptotic normality. 
We highlight that we are not restricted to the MLE
but we may consider as well 
the Generalized Method of Moment, pseudo- or quasi-likelihood
estimators, and other variants~(see \textit{e.g.},~\cite{Newey_2004}). 

Our procedure can also be used to study the lack of 
normality of the estimator under a fixed sample size $n$
when $G$ is Gaussian. 
The lack of normality may come from the fact that
$G_n$ has not reached the Gaussian regime and/or 
by the non-linearity of $H(n,\cdot)$. We develop
this in Example~\ref{exa:exponential} on
the exponential distribution and in Example~\ref{exa:binomial:2} on the binomial distribution.

Before introducing Theorem~\ref{th:random}, in the next section, we present some preliminary analytical results in which no randomness is involved.
However, to introduce some notations, we may refer to the desired probabilistic properties.

%%%%%%%%%%%%%%%%%%%%%%%%%%%%%%%%%%%%%%%%%%%%%%%%%%%%%%%%%%%%%%%%%%%%%%
\section{The asymptotic implicit function theorem}

Let us introduce the subsets $\Theta\subset\RR$ (the \emph{space of parameters}), 
$\II\subset\RR$, the \emph{space of rates}, 
and $\bS$, the \emph{space of scales}.

\begin{hypothesis}
\label{hyp:spaces}
The space $\Theta$ contains an open interval around~$0$, 
$(0,\epsilon)\subset\II$ for some $\epsilon>0$,
and there exists a function $\varphi:\bS \to \II$,
which we call the \emph{rate function}.
\end{hypothesis}

\begin{hypothesis}
    \label{hyp:beta}
We fix $p\in\Set{1,\dotsc,+\infty}$ and 
we consider 
$\bbeta:=\Set{\beta_k}_{k=0,1,\dotsc,p+1}\in\ZZ^{p+2}$.
\end{hypothesis}

\begin{notation}
    \label{not:delta}
For a parametric family of functions $\theta\mapsto F(s,\theta)$ which is of class $\cC^p(\Theta,\RR)$ 
for any $s\in\bS$ (or $F\in\cC^p(\Theta,\RR)^\bS$), 
we denote by $\dD^m$ the derivative of order $m$ with respect to~$\theta$.
We define 
\begin{equation*}
    \coef{F}_m(s):=\varphi(s)^{\beta_m} \dD^m F(s,0)\text{ for } m=0,\dotsc,p
    \text{ and }s\in\bS.
\end{equation*}
\end{notation}

The idea of the above notation is 
that the $\beta_m$ are such that $\varphi(s)^{\beta_m}\dD^m F(s,0)$ converges.

Let us introduce two notations, one related to $\bbeta$, and the other related
to multi-indices.

\begin{notation}
We define
\begin{align}
    \label{eq:betastar}
    \bbeta_\star&:=\max_{m=2,\dotsc,p+1 } \beta_m, 
    \\
    \label{eq:gamma}
\gamma_m&:=p+1+(\beta_m-\bbeta_\star)\wedge0
	\text{ for }m=1,\dotsc,p.
\end{align}
\end{notation}
%%

%%%
\begin{notation}
    In a systematic manner, we use sums overs multi-indices
     satisfying some constraints. Such multi-indices
     are denoted by $(k_1,\dotsc,k_m)$ and each $k_i$ is an integer. 
\end{notation}
%%%

We state a series expansion which approximates the roots of $F(s,\cdot)$.

%% THEOREM
%%
\begin{theorem}[Approximation of roots]
    \label{thm:1}
    Assume Hypotheses~\ref{hyp:spaces} and \ref{hyp:beta} hold and that~$\beta_1 \geq \bbeta_\star$.
    Let $F\in\cC^{p+1}(\Theta,\RR)^\bS$. Let $s\in\bS$ such that there exist a sequence 
	$\balpha(s)=\Set{\alpha_k(s)}_{k=0,\dotsc,p}$ 
	and a convex neighborhood $U(s)$ of $0$ in $\Theta$ satisfying: 
	\begin{enumerate}[thm]
	\item\label{thm:root:i} 
		 $\alpha_0(s)=0$  and
		\begin{multline}
			\label{eq:rs:s}
			\coef{F}_0(s) \varphi(s)^{-\beta_0}
			\\
			+\sum_{m=1}^p  \frac{\coef{F}_m(s)}{m!}
				\sum_{k_1+\dotsb + k_m < \gamma_m }
					\alpha_{k_1}(s)\dotsb\alpha_{k_m}(s)
    				\varphi(s)^{\sum_{n=1}^{m} k_n-\beta_m}=0.
    	\end{multline}

    \item\label{thm:root:ii} There exists 
	$\theta(s)\in U(s)$ such that $F(s,\theta(s))=0$.

    \item\label{thm:root:iii} It holds that
	\begin{gather*}
	\abs{\varphi(s)}^{\beta_1}\abs{\dD F(s,w)}\geq c(s)>0,
	\text{ for all }w\in U(s)\\
	\text{ and }
	\adjustlimits
		 \max_{m\in \{2,\ldots, p+1\}} \sup_{w\in U(s)} \abs{\dD^m F(s,w)}<+\infty.
	\end{gather*}

    \item\label{thm:root:iv}
    It holds that $\sum_{i=1}^p \abs{\alpha_i(s)}\cdot\abs{\varphi(s)}^i\in U(s)$.

    \end{enumerate}
Then 
\begin{equation*}
    \abs*{\theta(s)
	-
    \sum_{i=1}^p \alpha_{i}(s) \varphi(s)^i}
    \leq K(s) \max\Set*{
	\sum_{k=1}^p \abs{\alpha_k(s)},
    \Paren*{\sum_{k=1}^p \abs{\alpha_k(s)}}^{p+1}
} \abs{\varphi(s)}^{p+1},
\end{equation*}
with 
\begin{equation*}
    K(s)
    \leq \frac{\displaystyle \max_{\substack{m\in \{1,\ldots, p\} \\ 
	\beta_m +1\leq \bbeta_{\star}}} 
    \abs{\varphi(s)}^{\beta_1-\bbeta_{\star}} 
    \abs{\coef{F}_m(s)} + 
\adjustlimits\max_{m\in \{2,\ldots, p+1\}} \sup_{w\in U(s)} \abs{\varphi(s)}^{\beta_1} \abs{\dD^{m} F(s,w)} }{c(s)} <\infty.
\end{equation*}

\end{theorem}

The proof of Theorem~\ref{thm:1} is provided in Appendix~\ref{sec:asymopt-inversion}. 
We also consider the case of analytic functions, in which the expansion 
is infinite.

%% THEOREM
%%
\begin{theorem}[Approximation of roots, analytic case]
    \label{thm:1:analytic}
    Assume that Hypotheses~\ref{hyp:spaces} and \ref{hyp:beta} hold for $p=+\infty$.  
    Let $F:\bS\times\Theta\to\RR$ be such that $F(s,\cdot)$ is analytic.
    Let $s\in\bS$ such that there exist a sequence $\balpha(s)=\Set{\alpha_k(s)}_{k\geq 0}$ 
    and a convex neighborhood $U(s)$ of $0$ in $\Theta$ satisfying: 
    \begin{enumerate}[thm]
	\item\label{thm:root:i:analytic} 
	    $\alpha_0(s)=0$  and \eqref{eq:rs:s} holds with $p=+\infty$ and $\gamma_m=+\infty$
	    for $m\geq 1$.

    \item\label{thm:root:ii:analytic} There exists $\theta(s)\in U(s)$ such that $F(s,\theta(s))=0$.

    \item\label{thm:root:iii:analytic} It holds that 
	$
	    \abs{\varphi(s)}^{\beta_1}\abs{\dD F(s,w)}\geq c(s)>0
	$
 	 for all $w\in U(s)$.

\item\label{thm:root:iv:analytic}
    It holds that $\sum_{i=1}^{+\infty} \abs{\alpha_i(s)}\cdot\abs{\varphi(s)}^i\in U(s)$.
    \end{enumerate}
Then 
    \begin{equation*}
	\theta(s)= \sum_{i=1}^{+\infty} \alpha_{i}(s) \varphi(s)^i.
    \end{equation*}
\end{theorem}

In the next sections, in particular in Section~\ref{sec:alpha:construction} we define the notion of \emph{$(\bbeta,p)$-related sequence to a given sequence $\bdelta$} and we show that, if $\balpha$ is a $(\bbeta,p)$-related sequence to $\bdelta$, then it satisfies an equation analogous to~\eqref{eq:rs:s}:
		\begin{equation}
			\label{eq:rs:eval}
			\delta_0 \varphi(s)^{-\beta_0}
			+\sum_{m=1}^p  \delta_m
				\sum_{k_1+\dotsb + k_m < \gamma_m }
					\alpha_{k_1}\dotsb\alpha_{k_m}
    				\varphi(s)^{\sum_{n=1}^{m} k_n-\beta_m}=0.
    	\end{equation}
Note that~\eqref{eq:rs:s} corresponds to~\eqref{eq:rs:eval} with 
$\bdelta(s)= \Set*{{\coef{F}_k(s)}/{k!}}_{k=0,\ldots,p}$, 
but also with $\bdelta(s)=\Set*{-{\coef{F}_k(s)}/{(k! \coef{F}_1(s))}}_{k=0,\ldots,p}$ when $\coef{F}_1(s)\neq 0$.

Later, we focus on two sequences 
\begin{equation*}
    \bbeta^{\upflat}:=(1,2,2,\dotsc,2)
    \text{ and }
    \bbeta^{\zig}:=(1,2,1,2,\dotsc,1,2)
\end{equation*}
which are relevant for statistical applications. 
We refer to these choices of $\bbeta$ respectively as~\textquote{up-flat} and \textquote{zigzag} cases.

For the reader convenience, we provide here sequences $\balpha$ which are $(\bbeta,p)$-related sequence to a given sequence $\bdelta$ with $\delta_1=-1$ and $\bbeta=\bbeta^{\upflat}$ or $\bbeta^{\zig}$. 
For $\bbeta^{\upflat}$, let 
\begin{equation}
	\label{eq:upflat:eval}
	\alpha_1= \delta_0
	\quad \text{and} \quad 
    \alpha_{q} = \sum_{m=2}^{q} \delta_m \sum_{k_1+\dotsb+k_m=q}
    \alpha_{k_1}\dotsb \alpha_{k_m}
	\text{ for }2\leq q\leq p,
\end{equation}
and, for $\bbeta^{\zig}$, take
	$\alpha_1= \delta_0$, 
	$\alpha_{2q}=0$ for $2q\leq p$
    and
	for $2q+1\leq p$ 
    \begin{equation}
	\label{eq:zig:eval}
    \alpha_{2q+1}  
	 = \sum_{\substack{ m=2 \\ m \text{ even }}}^{2q}  
	 \delta_m  \sum_{k_1+\dotsb+k_m=2q} \alpha_{k_1} \cdots \alpha_{k_m}
	%\\	
	 + \sum_{\substack{ m=3 \\ m \text{ odd }}}^{2q+1} \delta_m
	 \sum_{k_1+\dotsb+k_m = 2q+1} \alpha_{k_1} \cdots \alpha_{k_{m}}.
\end{equation}
We refer to Lemmas~\ref{lem:rs},~\ref{lem:zigzag} and \ref{lem:upflat} in Section~\ref{sec:alpha:construction} for the proofs. %In the case \textquote{up-flat} the sequence can be given also by~Proposition~\ref{prop:inversion}.

For $\bbeta^{\upflat}$ with $p=5$ we have $\alpha_1={\delta}_0$ and
\begin{align*}
    \alpha_2  &= {\delta}_2 {\delta}_0^2,   &    \alpha_3&= ({\delta}_3 + 2 {\delta}_2^2) {\delta}_0^3,\\
    \alpha_4 &=({\delta}_4 + 5 {\delta}_2{\delta}_3 + 5 {\delta}_2^3) {\delta}_0^4, &   
    \alpha_5 &=({\delta}_5 + 6 {\delta}_2{\delta}_4 + 21 {\delta}_2^2{\delta}_3 + 14{\delta}_2^4 +3 {\delta}_3^2) {\delta}_0^5.
\end{align*}

For $\bbeta^{\zig}$ with $p=5$ we have $\alpha_1={\delta}_0$ and
\begin{align*}
    \alpha_2&= 0, & 
    \alpha_3&= {\delta}_3 {\delta}_0^3 +{\delta}_2 {\delta}_0^2,\\
    \alpha_4&=0, & 
    \alpha_5&= ({\delta}_5 + 3 {\delta}_3^2) {\delta}_0^5 + (5 {\delta}_2 {\delta}_3 + {\delta}_4) {\delta}_0^4 + 2 {\delta}_2^2 {\delta}_0^3.
\end{align*}
Note that $\alpha_{2k+1}$ is not homogeneous in ${\delta}_0^{2k+1}$
but contains terms in ${\delta}_0^{j}$  for $j \in \{ k+1, \ldots, 2 k+1\}$.

%%%%%%%%%%%%%%%%%%%%%%%%%%%%%%%%%%%%%%%%%%%%%%%%%%%%%%%%%%%%%%%%%%%%%%
%%%%%%%%%%%%%%%%%%%%%%%%%%%%%%%%%%%%%%%%%%%%%%%%%%%%%%%%%%%%%%%%%%%%%%
%%%%%%%%%%%%%%%%%%%%%%%%%%%%%%%%%%%%%%%%%%%%%%%%%%%%%%%%%%%%%%%%%%%%%%
%%%%%%%%%%%%%%%%%%%%%%%%%%%%%%%%%%%%%%%%%%%%%%%%%%%%%%%%%%%%%%%%%%%%%%
\section{The stochastic asymptotic implicit function theorem}

\label{sec:IFT}

In this section we consider the inference of a parameter $\theta_0\in\RR$ and 
we assume that~$\Theta$ contains an open interval around $\theta_0$.

To deal with convergence, we need some topological assumptions.
We consider that~$(\bS,\Yleft)$ is a directed set 
that is a set with a partial order $\Yleft$ 
such that for any $s,s'$ in $\bS$, there exists $r\in\bS$
such that $s\Yleft r$ and $s'\Yleft r$. 
Typically, we consider $\bS=\RR$ or $\bS=\NN$ (sample
size) yet this framework allows 
consider partitions and pairs (sample size and discretization step).

\begin{hypothesis} \label{hyp:rate}
    The rate function $\varphi:\bS\to\II\cap(0,1]$
    is monotone decreasing, that is $\varphi(s)\leq \varphi(s')$
    whenever $s'\Yleft s$.
\end{hypothesis}
Note that, since $(\bS,\Yleft)$ is a directed set, the rate function is a {\it net}.

\begin{notation} \label{not:deltaF}
Let $F\in\cC^{p}(\Theta,\RR)^\bS$ such that $\dD F(s,\theta_0)\neq 0$ then,
analogously to Notation~\ref{not:delta}, let for all $k\in \{0,\ldots,p\}$,
\begin{equation*}
	\coef{F}_k(s):=\varphi(s)^{\beta_k} \dD^k F(s,\theta_0).
\end{equation*}
Let $\bdelta^F$ be the family of $(\RR^{p+1})^{\bS}$
defined by 
\begin{equation*}
    \bdelta^F(s):=\Set*{ -\frac{1}{k !} \frac{ \coef{F}_k(s) }{ \coef{F}_1(s) }}_{k=0,\dotsc,p}.
\end{equation*}
\end{notation}

In the next result, our main result, we use the notion of $(\bbeta,p)$-related sequence to a given sequence $ \bdelta$. 
This notion has been briefly presented in the previous section, but it is introduced in full details in Section~\ref{sec:alpha:construction}.

%% PROPOSITION
\begin{proposition}
	\label{prop:limit:alpha}
	We consider a family of random functions $(F(s,\theta))_{(s,\theta)\in\bS\times \Theta}$
on a probability space $(\Omega,\cF,\PP)$
as well as $p$ as in Hypothesis~\ref{hyp:beta}.
	Assume that for all $s\in \bS$ and $\PP$-almost all $\omega \in \Omega$ the function $\theta\mapsto F(s,\theta)$
	is of class $\cC^{p}$ at $\theta_0$, and there exists a random vector $\bmu=\Set{\mu_k}_{k=0}^p$
	such that 
	 \begin{equation*}
	     (\coef{F}_0(s),\dotsc,\coef{F}_p(s)) \cvdist[s\to\infty] \bmu.
	 \end{equation*}
	For every $s\in \bS$ let $\balpha(s)$ the $(\bbeta,p)$-related sequence to $ \bdelta^F(s)$.
	Then
	\begin{enumerate}[thm]
	\item 
	    $ \bdelta^F(s) \cvdist[s\to\infty] \bdelta:=\Set*{- \frac{1}{m!}\frac{\mu_m}{\mu_1}}_{m=0,\dotsc,p}$.
	\item $\balpha(s) \cvdist[s\to\infty] \balpha$ where $\balpha$ is the $(\bbeta,p)$-related sequence to $\bdelta$.
    \end{enumerate}
\end{proposition}
%%%

%% DEFINITION
\begin{definition}[Boundedness in probability]
    A family of random variables $\Set{X_s}_{s\in\bS}$
    is said to be \emph{bounded in probability} if 
    for any $\epsilon>0$, there exist $K\geq 0$ and $s\in\bS$
    such that $\sup_{s\Yleft s'}\PP\Prb{\abs{X_{s'}}\geq K}\leq \epsilon$.
\end{definition}
%%

%% HYPOTHESIS
\begin{hypothesis}
    \label{hyp:random}
We consider a family of random functions $(F(s,\theta))_{(s,\theta)\in\bS\times \Theta}$
on a probability space $(\Omega,\cF,\PP)$
as well as $\bbeta\in\Set{\bbeta^{\upflat},\bbeta^{\zig}}$ and $p$ as in Hypothesis~\ref{hyp:beta}
and the rate function as in Hypothesis~\ref{hyp:rate}.
Assume that there exists an open neighborhood $U\subseteq \Theta$ of $\theta_0$ such that for all $s\in \bS$ and $\PP$-almost all $\omega \in \Omega$:
\begin{enumerate}[thm]
    \item\label{hyp:I} $\theta\mapsto F(s,\theta)$
	is of class $\cC^{p+1}$ in $U$.

     \item\label{hyp:II} There exists $\theta(s)\in U$
	such that $F(s,\theta(s))=0$.

    \item\label{hyp:III} $\dD F(s,\theta)\neq 0$
	for all $\theta\in U\setminus \{\theta_0\}$ and $\Set{1/\coef{F}_1(r)}_{r\in\bS}$ is bounded
	in probability.

    \item\label{hyp:IV} The family $\Set{N_r(F)}_{r\in\bS}$
	is bounded in probability with 
	\begin{equation*}
	    N_s(F):=\sup_{1\leq k\leq p+1}
	    \sup_{\theta\in U} \abs{\varphi(s)^{\beta_1}\dD^{k}F(s,\theta)}.
	\end{equation*}
\end{enumerate}
\end{hypothesis}

The following proposition is a direct consequences of Theorem~\ref{thm:1}.

%% THEOREM
\begin{theorem} 
    \label{th:random}
Assume Hypothesis~\ref{hyp:random} and assume that the net $\varphi$ in Hypothesis~\ref{hyp:rate} converges to~$0$. 
Assume that there exists a random vector $\bmu=\Set{\mu_k}_{k=0}^p$
	such that 
	\begin{equation*}
	    (\coef{F}_0(s),\dotsc,\coef{F}_p(s))
	    \cvdist[s\to\infty] \bmu,
	\end{equation*}
	and for every $s\in \bS$ let $\balpha(s)$ the $(\bbeta,p)$-related sequence to $\bdelta^F(s)$ (see~\eqref{eq:upflat:eval}-\eqref{eq:zig:eval}).
Let $\balpha$ and $\bdelta$ as in Proposition~\ref{prop:limit:alpha}.
We introduce
\begin{gather}
\pP{\ell}(\varphi(s); \balpha(s)) := \sum_{k=0}^\ell \alpha_k(s) \varphi(s)^k,
\\
\theta_\ell(s):=\theta_0+\pP{\ell}(\varphi(s),\balpha(s))
\text{  and  }
\theta_\ell(s,\infty) :=\theta_0+\pP{\ell}(\varphi(s),\balpha)
\end{gather}
for $\ell=0,\dotsc,p$.
Then,% using Notation~\ref{not:formal:3},
\begin{enumerate}[thm]
	\item $\pP{p}(\varphi(s); \balpha(s))$ converges to 0 in probability.
	\item $|\theta(s) - \theta_0 - \pP{p}(\varphi(s); \balpha(s)) | \varphi(s)^{-(p+1)}$ is bounded in probability.
	\item $\varphi(s)^{-(k+1)}  \Paren*{\pP{p}(\varphi(s); \balpha(s)) -\pP{k}(\varphi(s); \balpha(s))}  \cvdist[s\to\infty] \alpha_{k+1}$ and 
	\begin{equation*}    
	\varphi(s)^{-(k+1)}  (\theta(s) - \theta_k(s)) \cvdist[s\to\infty] \alpha_{k+1}
	\end{equation*}
	for $k \in \{0,\ldots, p-1\}$ (with $\pP{0} \equiv 0$).
		In particular if $k=0$ 
		\begin{equation*}
		    \varphi(s)^{-1}  (\theta(s) - \theta) \cvdist[s\to\infty] \alpha_{1}.
	    \end{equation*}
	\item $\theta(s) -\theta_p(s,\infty)$ 
		converges to 0 in probability with rate $\varphi(s)\left(\alpha_1(s)-\alpha_1\right)$.
\end{enumerate}
\end{theorem}
%%%%

The previous result exhibits two approximations of $\theta(s)$, \textit{i.e.},
\begin{equation*}
\theta_p(s):=\theta_0+\pP{p}(\varphi(s),\balpha(s))  \quad  \text{  and  }  \quad  \theta_p(s,\infty) :=\theta_0+\pP{p}(\varphi(s),\balpha).
\end{equation*}

In Section~\ref{sec:application} we consider applications in which the estimator $\theta(s)$ of $\theta_0 \in \Theta$ (assume~$\Theta$ is open) is a zero of 
$\theta \mapsto F(s,\theta):=g(\theta)T(s) - f(\theta)$
where $f,g\in \cC^p(\Theta)$
and~$T(s)$ is a sample statistic.
Indeed this is often the case for MLE estimators for independent observations with distributions belonging to an exponential family or for the method of moments. 

We aim at highlighting the following facts:
\begin{enumerate}[wide,label=(F\arabic*),topsep=-\parskip]
    \item When $F$ is the score, 
	the result above shows the importance of the hidden variable 
	$\alpha_1(s)=\delta^F_0(s)=-\coef{F}_0(s)/\coef{F}_1(s)$ to assess the quality of the MLE. 
	This variable $\alpha_1(s)$ depends on $F$ and its derivative. 
	The other variables $\alpha_k(s)$ are polynomial expressions
	in $\alpha_1(s)$ and $\delta^F_k(s)=-\coef{F}_k(s)/(k!\coef{F}_1(s))$ for $k\geq 1$.
	In the case \textquote{up-flat}, owing to the homogeneity
	in $\alpha_1(s)$ granted by~\eqref{eq:upflat:eval}, %Lemma~\ref{lem:upflat}, 
	\begin{equation*}
	    \theta(s)
	    =
	    \theta_0+\sum_{k=1}^{p} \xi_k(s)\alpha_1(s)^k\varphi(s)^k
	    +\grandO(\varphi(s)^{p+1}),
	\end{equation*}
	where the coefficients $\xi_k(s)$, $k\geq 2$ are polynomials expressions
	in the terms $\delta^F_\ell(s)$, %$=-\coef{F}_\ell(s)/\ell!\coef{F}_1(s)$, 
	$1\leq \ell\leq k$, 
	and $\xi_1(s)=1$. We shorten this as 
	\begin{equation}
	    \label{eq:mle:1}
	    \theta(s)=\theta_0+H\Paren*{s,\alpha_1(s)\varphi(s)}
	    +\grandO(\varphi(s)^{p+1}),
	\end{equation}
	where $H(s,\cdot)$ is a non-linear (polynomial or analytic) random function
	with $H(s,x)=x+\grandO(x^2)$.
	
	In many practical situations, the term $\coef{F}_0(s)/\coef{F}_1(s)=-\alpha_1(s)$ 
	fluctuates asymptotically as a (mixed) Gaussian random variable, 
	while the coefficients of the expansion of~$H(s,\cdot)$ converges in probability 
	to constants.

	A similar result holds for the case \textquote{zig-zag} ($\bbeta = \bbeta^{\zig}$), 
	excepted that the expression is no longer polynomial in $\alpha_1(s)\varphi(s)$,
	but still polynomial in $\alpha_1(s)$ and $\varphi(s)$ (granted by~\eqref{eq:zig:eval}). % Lemma~\ref{lem:zigzag}). 
	Some of the coefficients
	of $H$ converges in probability to (mixed) Gaussian distributions.

    \item Considering $\theta_p(s,\infty)$ instead of $\theta_p(s)$ 
	consists in replacing~$\alpha_1(s)$ and~$H(s,\cdot)$ by their limits.
	It gives then a \textquote{proxy} of the behavior of the MLE 
	which is useful when the convergence 
	of $\alpha_k(s)$'s toward their limits is faster than
	the one of $\theta_p(s)$. This happens with 
	the Skew Brownian motion \cite{lm22}, or with the exponential
	distribution as seen below in Example~\ref{exa:exponential}.

    \item In the examples in Section~\ref{sec:application} below, $\alpha_1(s)^2 \coef{F}_1(s)$ 
	(alternatively $\alpha_1(s)^2 \mu_1$ and $\coef{F}_0(s)^2/ \coef{F}_1(s)$ and $\coef{F}_0(s)^2/\mu_1$) 
	appears to be asymptotically \textquote{pivotal} 
	for $s$ large enough, in the sense that its distribution no longer depends 
	on the parameter. We could say it is the \emph{Wald's statistics}.
	Hence, we could change~\eqref{eq:mle:1}
	to express the MLE as a non-linear transform of one of these pivotal
	quantities. A further direction of research is to 
	use this to design confidence intervals or hypothesis 
	tests that takes the lack of normality into account.

    \item The upper error bound in Theorem~\ref{thm:1} (non-random version of Theorem~\ref{th:random}) also depends on $\alpha_1(s)$. 
	Therefore, when $\alpha_1(s)$
	deviates from $0$, the truncated expansions~$\theta_p(s)$ are not always accurate.

	\item 	
	    \label{fact:5}
	    In some situations, for $k\geq 1$, $\delta^F_k(s)$ %:=-\coef{F}_k(s)/k!\coef{F}_1(s)$
	    is both deterministic and constant in $s$ (see Section~\ref{sec:exp-fam:natural}).
	    In the case \textquote{up-flat} 
	    \begin{align*}
		\theta(s)&=\theta_0+H\Paren*{\delta^F_0(s)\varphi(s)} = %\theta_0+H\Paren*{-\frac{\coef{F}_0(s)}{\coef{F}_1(s)}\varphi(s)} = 
		\theta_0+H\Paren*{\alpha_1(s)\varphi(s)}
		\\
		\text{ and }
		\theta(s,\infty)&=\theta_0+H\Paren*{\alpha_1\varphi(s)}
	    \end{align*}
	    for a deterministic, polynomial (or analytic function $H$). 
	    The randomness in the MLE comes only from the one of $\coef{F}_0(s)$
	    that usually converges thanks to a CLT. 
		We have seen that $\varphi(s)^{-1}(\theta(s)-\theta_0)$ converges to $\alpha_1$, 
		and often in applications $\alpha_1$ has a Gaussian distribution $\cN(0,\sigma^2)$. 
		
	 Now, we quantify the distance of the distribution function of 
	 \begin{equation*}
		\varepsilon(s):=\varphi(s)^{-1}(\theta_p(s)-\theta_0)= \varphi(s)^{-1} H\Paren*{\alpha_1(s)\varphi(s)}
	    \end{equation*}
		from the one of its limiting distribution.
	    Let us consider the Kolmogorov-Smirnov distance:
	    For any two random variables $X_1,X_2$, let
		\begin{equation}
		    \label{eq:Kolm}
		\Delta(X_1;X_2):=\sup_{x\in \RR} \abs[\big]{\PP\Prb{ X_1 \leq x} - \PP\Prb{ X_2 \leq x} }.
		\end{equation}
		
	Whenever $H$ is invertible with inverse $H^{-1}$, we have
	\begin{equation*}
	\begin{split}
		\Delta(\varepsilon(s); \alpha_1 )
		& \leq \Delta(\varepsilon(s); \varphi(s)^{-1} H( \alpha_1 \varphi(s)) ) + \Delta(\varphi(s)^{-1} H( \alpha_1 \varphi(s)) ; \alpha_1)
		\\
		& = \Delta(\alpha_1(s); \alpha_1) +\Delta(\varphi(s)^{-1} H( \alpha_1 \varphi(s)) ; \alpha_1).
	\end{split}
	\end{equation*}
	The distance $\Delta(\varepsilon(s); \alpha_1 )$ is then given by the competition of two terms $\Delta(\alpha_1(s); \alpha_1)$ and 
	$\Delta(\varphi(s)^{-1} H( \alpha_1 \varphi(s)) ; \alpha_1)$.
	The quantity $\Delta(\alpha_1(s); \alpha_1) $ is the distance of $\varepsilon(s)$ to a possibly non linear transformation of $\alpha_1 \varphi(s)$.
	Indeed in the case we are considering~$H$ is deterministic and does not even depend on the sample size.
	The fact that $\Delta(\alpha_1(s); \alpha_1) \approx 0$ often follows from Berry-Esséen result or can be dealt with Stein's method,  
	and it rewrites:
	    for any $x\in\RR$, 
	    \begin{equation*}
		\PP\Prb*{\varepsilon(s)\leq x}
		\approx 
		    \PP\Prb*{\alpha_1 \leq \frac{H^{-1}(x\varphi(s))}{\varphi(s)}}.
	    \end{equation*}
	The term
	$\Delta(\varphi(s)^{-1} H( \alpha_1 \varphi(s)) ; \alpha_1)$, 
	with 
	\begin{equation*}
	    \frac{1}{\varphi(s)}H(\alpha_1\varphi(s))
	    = \alpha_1 + \delta_2 \alpha_1^2\varphi(s) + 
	(\delta_3+2\delta_2^2)\alpha_1^3\varphi(s)^2 
		+
	    \dotsb
	    ,
	\end{equation*}
	concerns the error introduced by the non-linearity.
	We show in Example~\ref{exa:exponential} that the effect of the latter error might be predominant.

	Note that we can also consider
	 \begin{align*}
		\theta(s)&=\theta_0+H\Paren*{\delta^F_0(s)\varphi(s)} %=\theta_0+H\Paren*{-\frac{\coef{F}_0(s)}{\coef{F}_1(s)}\varphi(s)}
		+\grandO(\varphi(s)^{p+1})
		\\
		\text{ and }
		\theta_p(s,\infty)&=\theta_0+H\Paren*{\alpha_1\varphi(s)},
	    \end{align*}
	    or dealing when $H$ depends on $s$ 
	    and is possibly random by using the limit of its coefficients as an approximation.

	\item Sometimes a suitable change of variable, as in the following Section~\ref{sec:change-var}, 
	    provides an expression of $\theta(s)$ instead of an approximation (see~Example~\ref{exa:binomial:2} below).

   \item There are also situations which exhibit some \emph{boundary layer},
	meaning that the coefficients $\alpha_k$ are themselves high
	and $\theta_p$ are not suitable approximation of the MLE. 
	This is the case for the binomial distribution with a parameter $\theta_0$
	close to~$0$ or~$1$. See Example~\ref{exa:binomial} or the article \cite{lm22}.

    \item For specific values of the parameter to estimate, sometimes several choices of $\bbeta$ are possible (note that $\beta^{\zig}_k \leq \beta^{\upflat}_k$), leading to different finite approximation of the estimator. We call this phenomena \emph{phase transition}. See Section~\ref{sec:symmetric} for an example or \cite{lm22}.
\end{enumerate}

%%%%%%%%%%%%%%%%%%%%%%%%%%%%%%%%%%%%%%%%%%%%%%%%%%%%%%%%%%%%%%%%%%%%%%
%%%%%%%%%%%%%%%%%%%%%%%%%%%%%%%%%%%%%%%%%%%%%%%%%%%%%%%%%%%%%%%%%%%%%%
%%%%%%%%%%%%%%%%%%%%%%%%%%%%%%%%%%%%%%%%%%%%%%%%%%%%%%%%%%%%%%%%%%%%%%
\section{Formal manipulation of expansions} \label{sec:alpha:construction}

In Theorems~\ref{thm:1} and \ref{thm:1:analytic}, 
the expansion of the root $\theta(s)$ to $F(s,\theta(s))=0$ is performed in term of 
$\varphi(s)$ with $s$-dependent coefficients. 
This implies a lack of uniqueness in the expansion.
To enforce uniqueness, we identify the coefficients
by transforming series to formal ones. 
Although natural, this choice is arbitrary
and \textquote{righteous}\footnote{As pointed out in \cite{Allaire2022}
regarding homogenization theory, a \textquote{criminal
path} leads to other ways to compute effective coefficients.}.

The main idea is to decouple the scale ($s$, playing the role of sample size) and the rate ($\varphi(s)$), considering the latter as a variable.
In Section~\ref{sec:dummy} we introduce the necessary operators to decouple and re-couple:
dummization and evaluation.
Next, in Section~\ref{sec:TFDB} we recall Faà di Bruno and Taylor formulae 
and show how dummization and evaluation behave in Taylor expansion of composition of functions.
In Section~\ref{sec:construction}, we propose a definition of sequences that satisfy Eq.~\eqref{eq:rs:s}. 
We call them $(\bbeta,p)$-related sequences.
We provide a different presentation in Section~\ref{sec:scramble}.
Next we show that there exists a unique $(\bbeta,p)$-related sequence 
when $\bbeta$ is $\bbeta^{\upflat}$ and~$\bbeta^{\zig}$ 
and we provide the proof of the recursive formulas~\eqref{eq:upflat:eval}-\eqref{eq:zig:eval} in Sections~\ref{sec:zigzag}-\ref{sec:upflat}.
In the \textquote{up-flat} case, we provide different characterizations of the $(\bbeta^{\upflat},p)$-related sequence 
(see Lemmas~\ref{lem:upflat} and~\ref{lem:composition:p})
that we propose also in the context of Theorem~\ref{thm:1} in Section~\ref{sec:seq:appl}.
To conclude, we show in Section~\ref{sec:seq:pert} how an additive perturbation acts on the sequence.

%%%%%%%%%%%%%%%%%%%%%%%%%%%%%%%%%%%%%%%%%%%%%%%%%%%%%%%%%%%%%%%%%%%%%%
\subsection{Dummization and evaluation} \label{sec:dummy}

For a ring $\KK$ and a family of indeterminates $\mathbf{X}:=(X_1,\dotsc,X_n)$ , 
$\KK\bra{X_1,\dotsc,X_n}$ (resp. $\KK\Paren{X_1,\dotsc,X_n}$)
denotes the ring of polynomials (resp. rational functions) with 
indeterminates $X_1,\dotsc,X_n$ and coefficients in $\KK$. 
Also, the set of formal power series is $\KK\dbra{X_1,\dotsc,X_n}$. 

To a finite, numerical sequence $\ba=(a_1,\dotsc,a_q)$, we consider 
the \textquote{dummization} operation which consists
in transforming each of the $a_i$'s to an \emph{indeterminate}
$a_i^{\dummy}$. Therefore $\RR\bra{\ba^{\dummy}}$ (resp. $\RR(\ba^{\dummy})$) is the 
ring of polynomials (resp. rational functions) whose indeterminates are the $a^{\dummy}_i$.

We consider another indeterminate $z$, that we will use later as a
placeholder for~$\varphi(s)^{\dummy}$ or for $x\in\RR$. 
We denote by $\KK\bra{z}$ the ring of polynomials whose coefficients
are in $\KK$. We are concerned with $\KK=\RR\Paren{\ba^{\dummy}}$. 
We also denote by $\KK\bra{z,z^{-1}}$ the Laurent series whose coefficients
are in~$\KK$. 

The \emph{evaluation} ---~denoted by~$\ce$~--- consists in transforming back $a_i^{\dummy}$ to $a_i$.
It is the dual operation of the dummization.
We extend $\ce$ as an algebra homomorphism
from~$\RR\dbra{\ba^{\dummy}}$ to~$\RR$. 
We furthermore extend $\ce$ on $\RR\dParen{\ba^{\dummy}}$
by setting $\ce(P/Q):=\ce(P)/\ce(Q)$ as well as to $\KK\dbra{z,z^{-1}}$
by setting $\ce(z^{k})=\varphi(s)^k$ for any $k\in\ZZ$.

\begin{notation}
    \label{not:formal:3}
Let $p\in \NN$. Given a sequence $\bdelta= \Set{\delta_k}_{k=0,\dotsc,p}\subset \KK$, 
we denote by~$\pP{p}(z;\bdelta)$ the polynomial
$\pP{p}(z;\bdelta):= \sum_{k=0}^p \delta_k z^k \in \KK\bra{z}$.
\end{notation}
\begin{notation}[Projection onto polynomials of degrees at most $p$]
    \label{not:formal:2}
Let $p\in \NN$.
We denote by $\fp_{p}$ is the projection of a 
power series $P\in \KK\dbra{z}$
onto the polynomials of degree at most $p$ by setting $z^q=0$ whenever $q>p$.
\end{notation}
Note that we can extend the latter notation to $p=\infty$ taking $\fp_{\infty}$ the identity
and~$\pP{\infty}(z;\bdelta)$ the series $\sum_{k=0}^\infty \delta_k z^k$ when it is well defined.

    The next lemma is immediate by expanding powers. 

\begin{lemma}
    \label{lem:composition}
	Let $p\in \NN$.
    Consider three families $\bdelta^{\dummy}$, 
	$\balpha^{\bullet}$, $\boldsymbol{\eta}^{\bullet}$ with components in a field $\KK$. 
	Assume that $\alpha^\bullet_0=0$.
	These families $\bdelta^{\dummy}$, $\balpha^\bullet$ and $\boldsymbol{\eta}^{\bullet}$ are linked by the relation 
    \begin{equation}
	\label{eq:fdb:1}
	\eta^{\bullet}_0=\delta_0^{\dummy}
	\text{ and }
	\eta^\bullet_k=\sum_{m=1}^k \delta_m^{\dummy}
	\sum_{k_1+\dotsb+k_m=k}
    \alpha^{\bullet}_{k_1}\dotsb\alpha^{\bullet}_{k_m}
    \in \RR(\bdelta^{\dummy})\text{ for all } k \in \{1,\ldots,p\}
    \end{equation}
	 if and only if 
	 \begin{equation*}
	    \pP{p}(z;\boldsymbol{\eta}^{\bullet}) = \fp_p \circ \pP{p}(\cdot;\bdelta^{\dummy}) \circ \pP{p}(z;\balpha^{\bullet}).
	    \end{equation*}
\end{lemma}
%%

%%%%%%%%%%%%%%%%%%%%%%%%%%%%%%%%%%%%%%%%%%%%%%%%%%%%%%%%%%%%%%%%%%%%%%
\subsection{Taylor and Faà di Bruno formulae} \label{sec:TFDB}

The Taylor formula is the main tool to go back
and forth between the realm of analysis and 
the one of algebra.

The Faà di Bruno formula is a way to compute the 
higher order derivatives of the composition 
of two functions \cite{zbMATH01956296}. It is a useful
tool to deal with power series \cite{abraham,johnston} and can be seen 
both from the points of view of analysis \cite{Leipnik2007}
or of algebra \cite{zbMATH06638469}.

\begin{proposition}[Faà di Bruno formula]
    \label{prop:fdb}
    Let $U,V$ be open subsets of $\RR$. 
    Let $f:U\to\RR$ and $g:V\to U$ be two functions of class $\cC^{p}$ for $p\geq 1$.
    Then $f\circ g$ is of class $\cC^{p}(V,\RR)$ 
	and~\footnote{We use a different convention than the usual ones on summing the indices
		written by regrouping the multi-index $(k_1,\dotsc,k_m)$
		by their respective values and counting the number of such partitions.}
    \begin{equation}
	\label{eq:fdb}
	\dD^\ell(f\circ g)(x)=\ell! \sum_{m=1}^\ell \frac{\dD^m f(g(x))}{m!}
	\sum_{k_1+\dotsb+k_m=\ell}
    \frac{\dD^{k_1} g(x)}{k_1!}\dotsb \frac{\dD^{k_m} g(x)}{k_m!},
    \end{equation}
    for $\ell =1,\ldots,p$ and $x\in V$.
\end{proposition}
\begin{definition}[Formal Taylor series]
    \label{def:formal:taylor}
    Let $f$ be a function of class $\cC^p(\RR,\RR)$
    for $p=1,\dotsc,+\infty$.
    We set $\delta_k:=\dD^k f(0)/k!$, $k=0,\dotsc,p$ and $\bdelta:=\Set{\delta_k}_{k=0,\dotsc,p}$. 
    The \emph{Taylor series operator} $\cT_p$ is defined as 
    \begin{equation*}
	\cT_p f(x):=\pP{p}(x;\bdelta)=\sum_{k=0}^p \delta_k x^k\in\RR\dbra{x}
    \end{equation*}
    that is, it gives the Taylor series of $f$ around $0$.
    The \emph{formal Taylor series} of~$f$ is
    \begin{equation}
	\label{eq:formal:5}
	{\cT_p^{\dummy} f}(z):=\pP{p}(z;\bdelta^{\dummy})=\sum_{k=0}^p \delta_k^{\dummy} z^k\in \KK\dbra{z}\text{ with } 
	    \KK=\RR\bra{\bdelta^{\dummy}} .
    \end{equation}
\end{definition}

If $\ce(z)=x$ for $f\in\cC^{p+1}(\RR,\RR)$,  
\begin{equation*}
    \ce({\cT_p^{\dummy} f}(z))
    =\cT_p f(x)=f(x)-R(x)
\end{equation*}
with $\abs{R(x)}\leq \abs{x}^{p+1} \sup_{y\in[-x,x]}\abs{\dD^{p+1} f(y)}/(p+1)!$
for any $x$ small enough.

A consequence of Lemma~\ref{lem:composition} and the Faà di Bruno formula
is the following one. 

\begin{lemma}
    \label{lem:taylor-compos}
    Let $p\in \NN$, let $f$, $g$ and $h$ be three functions of class $\cC^{p+1}$ such 
    that $g(0)=0$ and $\cT_p k(x)=\ce({\cT_p^{\dummy} k}(z))$ for $k=f,g,h$ with $\ce(z)=x$. 
    \\
	Then  
	${\displaystyle {\cT_p^{\dummy} h}= \fp_{p}( {\cT_p^{\dummy} f} \circ {\cT_p^{\dummy} g} )}$
	if and only if 
   $
	{\displaystyle \cT_p h
	=
	\cT_{p}(f\circ g)
	= \fp_p(\cT_p f\circ \cT_p g). }
    $
	\end{lemma}

We end this section by relating
the composition of analytic functions and their formal developments. 

\begin{notation}
    \label{not:analytic}
    For $r>0$, we denote by $\cA_r$ the class
    of analytic functions whose radius of convergence is
    at least $r$. We also set $B_r:=\Set{x\in\RR\given \abs{x}<r}$
    the open ball of radius $r$.
\end{notation}
\begin{hypothesis} 
    \label{hyp:infinity}
    The functions $f,g,h$ satisfy 
 $f\in\cA_{r_f}$, $g\in\cA_{r_g}$, 
  $h\in\cA_{r_h}$ with $g(B_{r_g})\subset B_{r_f}$, 
   so that $f\circ g\in\cA_{r_g}$.
\end{hypothesis}

Note that for an analytic function with a positive radius of convergence $r_f$, 
$f(x)=\cT_\infty f(x)$ whenever $\abs{x}<r_f$. 

\begin{lemma}
    \label{lem:taylor-compos-analytic}
    Let the same assumptions of Lemma~\ref{lem:taylor-compos}, and let $f,g,h$ as in Hypothesis~\ref{hyp:infinity}.
	Then  
	${\displaystyle {\cT_\infty^{\dummy} h}=  {\cT_\infty^{\dummy} f} \circ {\cT_\infty^{\dummy} g} }$
	if and only if 
   $
	{\displaystyle h
	=
	f\circ g.}
    $
\end{lemma}
\begin{remark} 
Let $p \in \NN$ and $\bdelta=\Set{\delta_k}_{k=0}^\infty \subset \RR$ such that
$\pP{\infty}(\cdot ;\bdelta) \in \cA_r$ for some $r>0$.
Then $\pP{p}(x;\bdelta)=\fp_p \pP{\infty}(x;\bdelta) = \cT_p \pP{\infty}(x;\bdelta)$ and $\pP{p}(x;\bdelta)= \cT_p \pP{p}(x;\bdelta)$.
In particular if $f \in \cA_r$ for some $r>0$, then
$\cT_p f= \fp_p\cT_\infty f = \cT_p\cT_\infty f$.
\end{remark}
%%

%%%%%%%%%%%%%%%%%%%%%%%%%%%%%%%%%%%%%%%%%%%%%%%%%%%%%%%%%%%%%%%%%%%%%%
%%%%%%%%%%%%%%%%%%%%%%%%%%%%%%%%%%%%%%%%%%%%%%%%%%%%%%%%%%%%%%%%%%%%%%
%%%%%%%%%%%%%%%%%%%%%%%%%%%%%%%%%%%%%%%%%%%%%%%%%%%%%%%%%%%%%%%%%%%%%%
\subsection{Constructing \texorpdfstring{$(\bbeta,p)$}{(β,p)}-related sequences}

\label{sec:construction}

The goal of this section is, for $s\in \bS$ fixed, to 
find a sequence $\balpha(s)$ satisfying~\eqref{eq:rs:s}.
The idea is to consider the coefficients~$\alpha_m(s)^{\bullet}$ and $\coef{F}_m(s)$ as dummy variables, then identify the coefficients $\alpha_m(s)^{\bullet}$ by the dummy version of Eq.~\eqref{eq:rs:s}, and finally transfer them back to their numerical values.

\begin{definition}[$(\bbeta,p)$-related sequences]
    \label{def:rs}
    Let $p$ and $\bbeta$ as in Hypothesis~\ref{hyp:beta}.
    Let $\bdelta=\Set{\delta_k}_{k=0,\dotsc,p}$
    and $\balpha=\Set{\alpha_k}_{k=0,\dotsc,p}$. 
    We say that $\balpha$ is \emph{$(\bbeta,p)$-related}
    to $\bdelta$ if  there exist~$\alpha^{\bullet}_k$, $k=0,\dotsc,p$
    which are rational functions in $\KK:=\RR\Paren{\bdelta^{\dummy}}$
    such that 
    \begin{enumerate}[thm]
	\item For $k=0,\dotsc,p$,  $\alpha_k=\ce(\alpha^{\bullet}_k)$
	    with $\alpha^{\bullet}_0=0$.
	\item In $\KK\dbra{z,z^{-1}}$, 
    \begin{equation}
	\label{eq:rs}
	\delta_0^{\dummy} z^{-\beta_0}
	+\sum_{m=1}^p  \delta_m^{\dummy}
	\sum_{k_1+\dotsb + k_m < \gamma_m }
	\alpha_{k_1}^{\bullet}\dotsb\alpha_{k_m}^{\bullet}
    z^{\sum_{n=1}^{m} k_n-\beta_m}=0,
    \end{equation}
    with $\gamma_m$ defined by \eqref{eq:gamma}, meaning 
    that all the coefficients in $\KK$ vanish for each
    of the powers of $z$.
    \end{enumerate}
\end{definition}
\begin{remark}
    \label{rem:shift}
    If $\balpha$ is a $(\bbeta,p)$-related sequence
    to $\bdelta$, then $\balpha$ it is also $(\bbeta',p)$-related
    sequence to $\bdelta$ for $\bbeta'=\Set{\beta_k+\gamma}_{k=0,\dotsc,p+1}$
    for any $\gamma\in\ZZ$.
\end{remark}
\begin{remark}
    \label{rem:scalar:multipl}
    If $\balpha$ is a $(\bbeta,p)$-related sequence
    to $\bdelta$, then $\balpha$ it is also $(\bbeta,p)$-related
    sequence to $\bdelta_\lambda$ for $\bdelta_\lambda=\Set{\lambda \delta_k}_{k=0,\dotsc,p}$
    for any $\lambda \in \RR$.
\end{remark}
\begin{example}[Non-existence/non-uniqueness]
    \label{ex:counterex1}
    We consider $p=2$, $\beta_0=0$, $\beta_1=1$, $\beta_2=2$ and $\beta_3\in\Set{1,2}$.
    Then $\gamma_m=1+\beta_m$ for $m=1,2$. Eq.~\eqref{eq:rs} becomes
    \begin{equation*}
        \delta^{\dummy}_0
	+
	\delta^{\dummy}_1\alpha^{\dummy}_1
	+
	\delta^{\dummy}_2(\alpha^{\dummy}_1)^2=0.
    \end{equation*}
    The possibility to find $\alpha^{\dummy}_1\in\RR$ depends on the relative
    values of $\bdelta$. Even if we are able to do so, 
    we are free to choose $\alpha^{\dummy}_2$ and the choice of $\balpha$
    which is $((0,1,2),2)$-related to $\bdelta$ is not unique.
\end{example}
\begin{example}[Non-existence/non-uniqueness]
    \label{ex:counterex2}
    We consider $p=2$, $\bbeta=(0,1,1,2)$. 
	If $\delta^{\dummy}_2 \neq 0$, then $\gamma_m=1+\beta_m$ for $m=1,2$. Eq.~\eqref{eq:rs} becomes
    \begin{equation*}
        \delta^{\dummy}_0
	+
	\delta^{\dummy}_1\alpha^{\dummy}_1
	+
	\delta^{\dummy}_2(\alpha^{\dummy}_1)^2z=0,
    \end{equation*}
    so that this could not be solved unless $\delta^{\dummy}_0=0$ and $\alpha^{\dummy}_1=0$.
    In the latter case, we are free to choose $\alpha^{\dummy}_2$.
\end{example}

The next lemma relates the definition above to~the assumption~\ref{thm:root:i} %(equation~\eqref{eq:rs:s}) 
in Theorem~\ref{thm:1}.
\begin{lemma}
    \label{lem:rs}
    Let $F\in\cC^p(\Theta,\RR)^{\bS}$ such that $\coef{F}_1\neq 0$, and fix $s\in\bS$.
    Let $\balpha(s)$ be a $(\bbeta,p)$-related sequence to
    $\bdelta^F(s)$. %with $\delta_k(s)=\delta^F_k(s) = - \coef{F}_k(s)/(k!\coef{F}_1(s))$, $k=0,\dotsc,p$. 
    Define the evaluation map so that $\ce(z)=\varphi(s)$. 
    Then $\balpha(s)$ satisfies~\eqref{eq:rs:s}. 
\end{lemma}
We study below two examples of $\bbeta$ for which $\balpha$ is
unique when using Definition~\ref{def:rs}. 
However, given a sequence $\balpha(s)$ satisfying~\eqref{eq:rs:s}, 
it is possible to change $\alpha_1(s)$ by a term of order $\varphi(s)$
so that \eqref{eq:rs:s} is still satisfied.
Therefore, there are an infinite number of solutions to \eqref{eq:rs:s}.
Our construction is a way to impose \emph{a unique choice} of~$\balpha(s)$
satisfying~\eqref{eq:rs:s} obtaining~\eqref{eq:upflat:eval} and~\eqref{eq:zig:eval}.
%%

%%%%%%%%%%%%

\subsection{Scrambling operator} \label{sec:scramble}

We provide here an alternative formulation of equation~\eqref{eq:rs} of Definition~\ref{def:rs}.

\begin{definition}[Scrambling operator]
    \label{def:scrambling}
    Given $\pP{p}(z;\bdelta)=\sum_{k=0}^p \delta_k z^k\in\KK\dbra{z}$ and~$\bbeta\in\RR^{p+1}$, 
    the \emph{scrambling operator} is
    \begin{equation*}
	\fs_{\bbeta}(\pP{p}(z;\bdelta)):=\sum_{k=0}^p \delta_k z^{k-\beta_k}\in \KK\dbra{z,z^{-1}}.
    \end{equation*}
\end{definition}

With Definition~\ref{def:scrambling}, Notation~\ref{not:formal:3} and Notation~\ref{not:formal:2}, 
\eqref{eq:rs} is equivalent to 
\begin{equation} \label{eq:rs:scramble}
    \fs_{\bbeta}\circ \fp_p \circ \pP{p}(\cdot;\bdelta^{\dummy}) \circ \pP{p}(z;\balpha^{\bullet})=0
    \text{ with }
\ \alpha^{\bullet}_0=0.
\end{equation}

Under the assumptions of Lemma~\ref{lem:rs}, \eqref{eq:rs} (and so \eqref{eq:rs:scramble}) rewrites
\begin{equation} \label{eq:formal:6}
\fs_{\bbeta}\circ \fp_p\circ {\cT_p^{\dummy} F} \Paren*{s,\pP{p}(z;\balpha(s)^{\bullet})}=0  \text{ in }\RR\Paren{\bdelta^{\dummy}(s)}\dbra{z,z^{-1}}.
\end{equation}

We now study $(\bbeta,p)$-related sequences for
some particular choices of $\bbeta$, which are
fitted for statistical applications.

%%%%%%%%%%%%%%%%%%%%%%%%%%%%%%%%%%%%%%%%%%%%%%%%%%%%%%%%%%%%%%%%%%%%%%
%%%%%%%%%%%%%%%%%%%%%%%%%%%%%%%%%%%%%%%%%%%%%%%%%%%%%%%%%%%%%%%%%%%%%%
%%%%%%%%%%%%%%%%%%%%%%%%%%%%%%%%%%%%%%%%%%%%%%%%%%%%%%%%%%%%%%%%%%%%%%
\subsection{The case \textquote{zig-zag}} \label{sec:zigzag}

Our first case of interest is $\bbeta^{\zig}=(1,2,1,2,\dotsc)$.
Its statistical usefulness is justified below in Section~\ref{sec:symmetric}. 
Eq.~\eqref{eq:rs} becomes  
\begin{multline*}
    \delta^{\dummy}_0  z 
    +\sum_{\substack{m\geq 0\\2m+1\leq p}}
    \delta^{\dummy}_{2m+1} 
    \sum_{k_1+\dotsb+k_{2m+1}\leq p} \alpha^{\bullet}_{k_1}\dotsb \alpha^{\bullet}_{k_{2m+1}}
    z^{k_1+\dotsb+k_{2m+1}}
    \\
    +
    \sum_{\substack{m\geq 1\\2m\leq p}}
    \delta^{\dummy}_{2m} 
    \sum_{k_1+\dotsb+k_{2m}<p} \alpha^{\bullet}_{k_1}\dotsb \alpha^{\bullet}_{k_{2m}}
    z^{k_1+\dotsb+k_{2m}+1}
    =0.
\end{multline*}

And it can be rewritten using the scrambling operator, which in this case is
\begin{equation*}
    \fs_{\bbeta^{\zig}}\Paren*{\sum_{k=0}^{+\infty} \delta^{\dummy}_k z^k}
    =
    \frac{1}{z^2}
    \Paren*{
    \sum_{k=0}^{+\infty} \delta^{\dummy}_{2k+1}z^{2k+1}
    +z\sum_{k=0}^{+\infty} \delta^{\dummy}_{2k}z^{2k}
}.
\end{equation*}

\begin{lemma}
    \label{lem:zigzag}
    %Assume that $\delta_1\neq0$. 
	Assume that $\delta_1=-1$.
    Then $\balpha$ is $(\bbeta^{\zig},p)$-related to $\bdelta$
    if and only if 
	$\alpha^{\bullet}_1= \delta^{\dummy}_0$, 
	$\alpha^{\bullet}_{2q}=0$ for $2q\leq p$
    and 
    \begin{multline*}
    \alpha^{\bullet}_{2q+1}  
	 = \sum_{\substack{ m=2 \\ m \text{ even }}}^{2q}  
	 \delta^{\dummy}_m \sum_{k_1+\dotsb+k_m=2q} \alpha^{\bullet}_{k_1} \cdots \alpha^{\bullet}_{k_m} 
		 + \sum_{\substack{ m=3 \\ m \text{ odd }}}^{2q+1} \delta^{\dummy}_m  
	 \sum_{k_1+\dotsb+k_m = 2q+1} \alpha^{\bullet}_{k_1} \cdots \alpha^{\bullet}_{k_{m}}
    \end{multline*}
	for $2q+1\leq p$.
%	$\alpha^{\bullet}_1=-\delta^{\dummy}_0/\delta^{\dummy}_1$, 
%	$\alpha^{\bullet}_{2q}=0$ for $2q\leq p$
%    and 
%    \begin{multline*}
%    \alpha^{\bullet}_{2q+1}  
%	 = - \sum_{\substack{ m=2 \\ m \text{ even }}}^{2q}  
%	 \frac{\delta^{\dummy}_m}{\delta^{\dummy}_1}\sum_{k_1+\dotsb+k_m=2q} \alpha^{\bullet}_{k_1} \cdots \alpha^{\bullet}_{k_m} 
%	\\	
%	 - \sum_{\substack{ m=3 \\ m \text{ odd }}}^{2q+1} \frac{\delta^{\dummy}_m}{\delta^{\dummy}_1}  
%	 \sum_{k_1+\dotsb+k_m = 2q+1} \alpha^{\bullet}_{k_1} \cdots \alpha^{\bullet}_{k_{m}}
%    \text{ for }2q+1\leq p.
%    \end{multline*}
    In particular, there exists one and only one sequence
    $\balpha$ which is $(\bbeta^{\zig},p)$-related to~$\bdelta$.
\end{lemma}

Lemma~\ref{lem:rs} and Lemma~\ref{lem:zigzag} prove~\eqref{eq:zig:eval}.

%%%%%
\begin{remark} \label{rem:infty-p-related:zig}
%Let $p\neq \infty$.
If $\balpha$ is $(\bbeta^{\zig},\infty)$-related to $\bdelta$,
then $\Set{\alpha_k}_{k=0}^p$ is $(\bbeta^{\zig},p)$-related to $\Set{\delta_k}_{k=0}^p$.
\end{remark}
%%%%%%%%%

%%
\begin{lemma}
We consider a sequence $\balpha$ which is $(\bbeta^{\zig},+\infty)$-related to $\bdelta$. 
Let $p\in \{1,2,\ldots, \infty\}$. %, $\balpha$ be $(\bbeta^{\upflat},p)$-related to $\bdelta$, 
We set 
\begin{equation*}
    f_{\mathrm{odd}}(x):=\sum_{\substack{k\geq 0\\2k+1\leq p}} \delta_{2k+1} x^{2k+1}
    \quad \text{ and } \quad 
    f_{\mathrm{even}}(x):=\sum_{\substack{k\geq 1\\2m\leq p}} \delta_{2k} x^{2k}
\end{equation*}
and $g(x):=\pP{p}(x;\balpha)$.
Assume in addition %Hypothesis~\ref{hyp:infinity} 
that $f_{\mathrm{odd}}, f_{\mathrm{even}}\in \cA_{r_f}$
and $g\in \cA_{r_g}$ with $g(B_{r_g})\subset B_{r_f}$ (necessary if $p=\infty$).
Then 
\begin{equation}
    \label{eq:zig:2}
    \cT_p\left( f_{\mathrm{odd}}(g(x))+x f_{\mathrm{even}}(g(x)) \right)= -\delta_0 x\text{ for }x\in B_{r_f}.
\end{equation}
Besides, if $\delta_0 \delta_1\neq 0$, then there exists an analytic function $h$
that solves \eqref{eq:zig:2} (with~$h$ instead of $g$) and therefore $\alpha_k:=\dD^k h(0)/k!$ for any $k\geq 1$.
\end{lemma}
\begin{proof}
We restrict to the case $p=\infty$, see Remark~\ref{rem:infty-p-related:zig}.
Let us introduce a new parameter $\nu$ and we define 
\begin{equation*}
    \widecheck{f}(x,\nu):=f_{\mathrm{odd}}(x)+\nu f_{\mathrm{even}}(x)
\end{equation*}
and the analogous $\widecheck{f}^{\dummy}(z,\nu)$ with $\bdelta^{\dummy}$ 
and $g^{\dummy}(z):=\pP{\infty}(z;\alpha^{\bullet})$.
As in Lemma~\ref{lem:composition}, 
if $\balpha^{\bullet},\bgamma^{\bullet}\in\RR\Paren{\bdelta^{\dummy}}$
and $\bdelta^{\dummy}$ are linked by the relation 
\begin{multline}
    \label{eq:zig:1}
    \sum_{m\geq 0} \delta^{\dummy}_{2m+1} \sum_{k_1,\dotsc,k_{2m+1} \geq 0} 
    \alpha^{\bullet}_{k_1}\dotsb\alpha^{\bullet}_{k_{2m+1}}z^{k_1+\dotsb+k_{2m+1}}
    \\
    +\nu^{\dummy}\sum_{m\geq 1} \delta^{\dummy}_{2m} \sum_{k_1,\dotsc,k_{2m}\geq 0} 
    \alpha^{\bullet}_{k_1}\dotsb\alpha^{\bullet}_{k_{2m}}z^{k_1+\dotsb+k_{2m}}
    \\
    =\sum_{m\geq 0}\gamma_{2m+1}^\bullet z^{2m+1}
    +\nu^{\dummy}\sum_{m\geq 1}\gamma_{2m}^\bullet z^{2m}
\end{multline}
in $\RR\Paren{\delta^{\dummy}}\dbra{z,\nu^{\dummy}}$ for some indeterminate $\nu^{\dummy}$,
then 
\begin{equation*}
    \widecheck{f}^{\dummy}(g^{\dummy}(x),\nu^{\dummy}) 
	=
    \sum_{m\geq 0}\gamma_{2m+1}^{\bullet} x^{2m+1}
    +\nu^{\dummy}\sum_{m\geq 1}\gamma_{2m}^{\bullet} x^{2m}.
\end{equation*}
If $\balpha$ is $(\bbeta^{\zig},+\infty)$-related to $\bdelta$, 
then \eqref{eq:zig:1} holds for $\nu^{\dummy}=z$, $\gamma_k=0$ for $k=0,2,3,\dotsc$
and $\gamma_1=-\delta_0$ with $\gamma_k:=\ce(\gamma^{\bullet})$. This proves~\eqref{eq:zig:2}. 

We set $q(y,x)=\delta_0 x+\widecheck{f}(y,x)$.
If $\delta_1\neq 0$, then $D_y q(y,0)\neq 0$. 
From the analytic implicit function theorem (see \textit{e.g.},~\cite{abraham}),
there exists a unique analytic function $h$
defined in a neighborhood of $0$ 
such that $q(h(x),x)=0$, meaning
that $h$ solves 
$f_{\mathrm{odd}}(h(x))+x f_{\mathrm{even}}(h(x)) = - \delta_0 x$ and so \eqref{eq:zig:2}.
By uniqueness in a neighborhood of 0, $h=g$.
\end{proof}
%%

%%%%%%%%%%%%%%%%%%%%%%%%%%%%%%%%%%%%%%%%%%%%%%%%%%%%%%%%%%%%%%%%%%%%%%
%%%%%%%%%%%%%%%%%%%%%%%%%%%%%%%%%%%%%%%%%%%%%%%%%%%%%%%%%%%%%%%%%%%%%%
%%%%%%%%%%%%%%%%%%%%%%%%%%%%%%%%%%%%%%%%%%%%%%%%%%%%%%%%%%%%%%%%%%%%%%
\subsection{The case \textquote{up-flat}} \label{sec:upflat}

We consider the case $\bbeta^{\upflat}=(1,2,\dots,2)$, which 
we call the case \textquote{up-flat}.
With Remark~\ref{rem:shift}, \eqref{eq:rs} in the definition of $(\bbeta^{\upflat},p)$-related sequences  (Definition~\ref{def:rs}) can be replaced by
\begin{equation}
    \label{eq:upflat:1}
    \delta_0^{\dummy}  z 
    +\sum_{m=1}^ p \delta_m^{\dummy} \sum_{k_1+\dotsb + k_m\leq p} \alpha^\bullet_{k_1}\dotsb \alpha^\bullet_{k_m}z^{k_1+\dotsb+k_m}=0.
\end{equation}

We now give an explicit, recursively computable expressions for the
$\alpha_k$'s from the~$\delta_k$.  Other ways to express the $\alpha_k$'s are
possible, for example using Bell polynomials.  See also
\cite{zbMATH01891065,Brent1978} for considerations on the efficiency of
numerical implementations.

\begin{lemma}
    \label{lem:upflat}
    Assume that $\delta_1\neq0$.
    Then $\balpha$ is $(\bbeta^{\upflat},p)$-related to $\bdelta$
    if and only if 
   \begin{equation}
	\label{eq:upflat:4}
	\alpha^\bullet_1=\frac{-\delta^{\dummy}_0}{\delta^{\dummy}_1}
	\text{ and }
    \alpha^\bullet_{q} =-\sum_{m=2}^{q} \frac{\delta^{\dummy}_m}{\delta^{\dummy}_1} \sum_{k_1+\dotsb+k_m=q}
    \alpha^\bullet_{k_1}\dotsb \alpha^\bullet_{k_m}
    \text{ for }2\leq q\leq p.
    \end{equation}
    In particular, there exists one and only one sequence
    $\balpha$ which is $(\bbeta^{\upflat},p)$-related to~$\bdelta$.
    Besides, each of the $\alpha_q^\bullet$ is a polynomial 
    expression in $\delta_m^{\dummy}$, $m=0,\dotsc,q$
    whose term in $\delta_0^{\dummy}$ has degree $q$.
\end{lemma}
\begin{proof} We rewrite \eqref{eq:upflat:1} as 
\begin{equation}
    \label{eq:upflat:1bis}
    \delta^{\dummy}_0  z 
    +\sum_{q=1}^p \sum_{m=1}^ p \delta^{\dummy}_m 
    \Paren*{\sum_{k_1+\dotsb + k_m =q}
    \alpha^\bullet_{k_1}\dotsb \alpha^\bullet_{k_m}}
    z^{q}=0.
\end{equation}
We identify $\alpha^\bullet_k$'s according to the power of $z$. 
In particular, $\delta^{\dummy}_0+\delta^{\dummy}_1 \alpha^\bullet_1=0$ 
so that $\alpha^\bullet_1=-\delta^{\dummy}_0/\delta^{\dummy}_1$.
Since $\alpha^{\dummy}_0=0$, for $2\leq q\leq p$, 
\begin{equation*}
    0=\sum_{m=1}^p \delta^{\dummy}_m\sum_{k_1+\dotsb + k_m =q} \alpha^\bullet_{k_1}\dotsb \alpha^\bullet_{k_m}
    =
    \delta^{\dummy}_1 \alpha^\bullet_q+
    \sum_{m=2}^q\delta^{\dummy}_m 
    \sum_{\substack{1\leq k_1,\dotsc,k_m<q\\\
    k_1+\dotsb + k_m =q}}
    \alpha^\bullet_{k_1}\dotsb \alpha^\bullet_{k_m}.
\end{equation*}
This proves that a $(\bbeta^{\upflat},p)$-related sequence may be
constructed for any $\bdelta$ with $\delta_1\neq 0$.
The converse is proved by the same argument of identifying the coefficients
in the power of $z$.

The form of  the expression and the degree of the term $\delta_0^{\dummy}$
is proved by a direct induction.
\end{proof}

Lemma~\ref{lem:rs} and Lemma~\ref{lem:upflat} prove~\eqref{eq:upflat:eval}.

We now consider an alternative way to compute the coefficients in terms of the inverse of an analytic function.
Let us consider the scrambling operator and~\eqref{eq:rs:scramble}, which was an alternative expression to~\eqref{eq:rs}, in this case becomes the following alternative expression to~\eqref{eq:upflat:1}:
\begin{equation*}
0=\fs_{\bbeta^{\upflat}-2}\circ \fp_p \circ \pP{p}(\cdot;\bdelta^{\dummy}) \circ \pP{p}(z;\balpha^{\bullet}) =  \fp_p \circ \pP{p}(\cdot;\bdelta^{\dummy}) \circ \pP{p}(z;\balpha^{\bullet}) - \delta_0^{\dummy} + \delta_0^{\dummy} z.
\end{equation*}
Therefore,  \eqref{eq:rs} in Definition~\ref{def:rs} (represented here by~\eqref{eq:upflat:1}) can be replaced by
\begin{equation}
    \label{eq:upflat:2} 
	\fp_p \circ f^{\dummy} \circ \pP{p}(z;\balpha^{\bullet})= -\delta_0^{\dummy}z
\end{equation}
with $f^{\dummy}(z):=\pP{p}(z;\bdelta^{\dummy}) - \delta_0^{\dummy}=\sum_{k=1}^p \delta_k^{\dummy} z^k$.

\begin{remark}
To prove the equivalence of~\eqref{eq:upflat:1} and~\eqref{eq:upflat:2} we could also use Lemma~\ref{lem:composition}, if $p\in \NN$.
Indeed~\eqref{eq:upflat:1} is equivalent to~\eqref{eq:fdb:1} in Lemma~\ref{lem:composition} 
with $\eta_0=\delta_0$, $\eta_1=-\delta_0$  and $\eta_k=0$ for $k\geq 2$
(let us recall that $\alpha_0=0$ by definition).
\end{remark}

We then see that to get $\balpha$, we compute the compositional inverse of a
power series. Hence, our problem is related to the Lagrange
inversion formula~\cite{gessel}.   Another consequence is that one
may look at other ways than series expansions to compute approximations of the
root, for example using continued fractions or close form expressions.

The next lemma, Lemma~\ref{lem:composition:p}, ensures that identifying the inverse of a function 
constructed as a series/polynomial of coefficients $\delta_k$'s
yields the coefficients $\alpha_k$'s.

\begin{lemma} \label{lem:composition:p}
Let $p\in \{1,2,\ldots, \infty\}$, $\balpha$ be $(\bbeta^{\upflat},p)$-related to $\bdelta$, 
and
$f(x):=\sum_{k=1}^p \delta_k x^k=\pP{p}(x;\bdelta)-\delta_0$ and $g(x):=\pP{p}(x;\balpha)$.
If $p=\infty$, assume in addition Hypothesis~\ref{hyp:infinity}.
Recall that $g(0)=0$.
Then
\begin{equation*}
\cT_p(f \circ g (x)) = \fp_p ( f \circ  g (x)) =- \delta_0 x. 
\end{equation*}
If in addition $\delta_0\neq 0$ and $\delta_1 \neq 0$ then $f$ is invertible on an open neighborhood of $0$ and
\begin{equation*}
\alpha_k= (-\delta_0)^k \dD^k f^{-1}(0)/k! \quad \text{for any} \quad k\geq 1. 
\end{equation*}
\end{lemma}

\begin{proof}
The first part is a direct consequence of the combination of dummization-evaluation, Eq.~\eqref{eq:upflat:2}
(see Lemmas~\ref{lem:taylor-compos}-\ref{lem:taylor-compos-analytic}).
The second part of the statement follows from 
the Inverse Function Theorem for analytic functions \cite{abraham} which
one may applies when $\delta_1=\dD f(0)\neq 0$.
The proof is completed, by uniqueness shown in Lemma~\ref{lem:upflat}.
\end{proof}

The next corollary follows from the fact that $\cT_p ( f \circ  g) = \fp_p \left( \cT_p f \circ \cT_p g \right)$.

%% COROLLARY
\begin{corollary} \label{rem:composition:p}
Let $p\in \{1,2,\ldots, \infty\}$, $\balpha$ be $(\bbeta^{\upflat},p)$-related to $\bdelta$, 
and let $f,g$ functions of class $\cC^p$ (if $p=\infty$, assume in addition Hypothesis~\ref{hyp:infinity}) 
such that 
$\cT_p f (x) = \sum_{k=1}^p \delta_k x^k =\pP{p}(x;\bdelta)-\delta_0$ and $\cT_p g(x)=\pP{p}(x;\balpha)$.
Then
\begin{equation*}
\cT_p ( f \circ  g (x)) =- \delta_0 x. 
\end{equation*}
\end{corollary}
%%%%%%%%%

%% REMARK
\begin{remark} \label{rem:infty-p-related}
Let $p\neq \infty$.
If $\balpha$ is $(\bbeta^{\upflat},\infty)$-related to $\bdelta$,
then $\Set{\alpha_k}_{k=0}^p$ is $(\bbeta^{\upflat},p)$-related to $\Set{\delta_k}_{k=0}^p$.
\end{remark}
%%%%%%%%%

The above remark ensures that we can take always $p=+\infty$. For instance taking~$\delta_k=0$
for $k$ above a given threshold, say $p$, and we associate with $\bdelta$ 
an infinite sequence $\balpha$.

%% EXAMPLE
\begin{example}
    \label{ex:1}
    Let $\delta_k=(‐1)^k/\theta^{k-1}$ (resp.~$\delta_k=1/\theta^{k-1}$) for $k\geq 1$ and some $\theta\in\RR$. 
	Then 
    for $\abs{x}<\abs{\theta}$, 
    \begin{equation*}
	f(x):=\sum_{k=1}^\infty \delta_k x^k=-\frac{\theta x}{\theta+x} \quad (\text{resp.}~f(x)=\frac{\theta x}{\theta-x}).
    \end{equation*}
    Since $f^{-1}(x)=f(x)$ (resp. $f^{-1}(x)=-f(-x)= \frac{\theta x}{\theta+x}$) then
	$ \alpha_k=\frac{\delta_0^k}{\theta^{k-1} }$ (resp.~$\alpha_k=-\frac{\delta_0^k}{\theta^{k-1}  }$ ) for $k\geq 1$.
\end{example}
%%

%% EXAMPLE
\begin{example}
\label{ex:2}
Let $\psi$ be an analytic function around $0$ with $\psi'(0)\neq 0$
(so that $\psi$ is invertible around $0$) and set
\begin{equation*}
    \delta_k:=-\frac{1}{k!}\frac{\dD^k\psi(0)}{\psi'(0)}
    \text{ for }k\geq 1
\end{equation*}
 as well as 
\begin{equation*}
    f(z):=\sum_{k\geq 1} \delta_k z^k=\sum_{k\geq 1}\frac{-1}{k!}\frac{\dD^k\psi(0)}{\psi'(0)}z^k
    =\frac{\psi(0)-\psi(z)}{\psi'(0)}.
\end{equation*}
Then 
$
f^{-1}(z)=\psi^{-1}(\psi(0)-z\psi'(0))
$. 
By Faà di Bruno formula~\eqref{eq:fdb}
$
    \dD^k f^{-1}(z)=\dD^k \psi^{-1}(\psi(0)- z\psi'(0))\psi'(0)^k
$
and 
\begin{equation*}
    \alpha_k=  \frac{(- \psi'(0) \delta_0)^k \dD^k \psi^{-1}(\psi(0))}{k!}
    =  \frac{\psi(0)^k \dD^k \psi^{-1}(\psi(0))}{k!}.
\end{equation*}

\end{example}
%%
%%%

We now give a non trivial example.

\begin{example}
    Consider $\theta\not\in\Set{0,1}$.
    Let us assume that $\delta_1=-1$ and $\delta_k=(-\theta)^k$, $k\geq 2$.
    Then for $\abs{z}<1/\abs{\theta}$,
    \begin{equation*}
	f(z)
	=- z + \sum_{k\geq 2}(-\theta)^kz^k
	=
	\frac{1}{1+\theta z}-1+(\theta-1) z.
    \end{equation*}
    To compute the inverse $f^{-1}$, we are looking for the solution $h(z)$ to the quadratic equation
    \begin{equation*}
	\frac{1}{1+\theta h(z)}+(\theta-1) h(z)=1- z,\ \text{ for any $z$ small enough}. 
    \end{equation*}
	We assume $\theta \in (0,1)$.
	The possible solutions are
	\begin{equation*}
		h_\pm(z) =\frac{1- \theta z \pm \sqrt{ \theta^2 z^2 - 4 \theta^2 z + 2 \theta z+1}}{2\theta(\theta-1)},
	\end{equation*}
	where the quantity in the square root is always positive.
    But $f(0)=0$ hence $f^{-1}=h_-$.
    One could then check that $h_-$ has the same series expansion in $z$
    as the one in Lemma~\ref{lem:upflat}. It is also possible to express
     $h_-$ using a continued fraction.
\end{example}
\begin{corollary}
Let $r>0$ and let $\Set{f_n}_n$ be a sequence of functions in $\cA_r$
and $f\in \cA_r$ be such that $f_n$ converge uniformly to $f$. 
Then~$\delta^n_k:=\dD^k f_n(0)/k!$ converges to $\delta_k$ 
for any $k\geq 1$. Besides, if $\delta_1\neq 0$, then 
for $n$ large enough, the corresponding sequence~$\balpha^n$
which is $(\bbeta^{\upflat},+\infty)$-related to $\bdelta^n=\Set{\delta_k^n}_{k\geq 0}$
converges to $\balpha$, the sequence which
is~$(\bbeta^{\upflat},+\infty)$-related to $\bdelta=\Set{\delta_k}_{k\geq 0}$.
\end{corollary}
%%

%%%%%%%%%%%%%%%%%%%%%%%%%%%%%%%%%%%%%%%%%%%%%%%%%%%%%%%%%%%%%%%%%%%%%%

\subsection{\texorpdfstring{$(\bbeta,p)$}{(β,p)}-related sequences in the context of Theorem~\ref{thm:1}} 

\label{sec:seq:appl}

We have already seen in Lemma~\ref{lem:rs} 
how to derive a sequence $\balpha(s)$ satisfying~\eqref{eq:rs:s}:
by finding a $(\bbeta,p)$-related sequence to $\bdelta^F(s)$. % with $\delta_k(s)=\coef{F}_k(s)/k!$, $k=0,\dotsc,p$. 
Hence, when $\bbeta=\bbeta^{\upflat}$ (resp.~$\bbeta=\bbeta^{\zig}$), 
Lemma~\ref{lem:upflat} (resp.~Lemma~\ref{lem:zigzag}) provides a recursive
expression of the desired sequence that we have given in~\eqref{eq:upflat:eval} (resp.~\eqref{eq:zig:eval}).
We now apply Lemma~\ref{lem:composition:p} to provide a different expression in the case $\bbeta=\bbeta^{\upflat}$ in this context.

%% PROPOSITION
\begin{proposition}
    \label{prop:inversion}
    Let $F$ be a function of class $\cC^{p+1}(\Theta,\RR)^\bS$.
    Take $\theta_0\in\Theta$ 
	and consider Notation~\ref{not:deltaF}. 
 %   $\coef{F}_k:=\varphi(s)^2\dD^k F(s,\theta_0)$ for $k\geq 1$
 %  and $\coef{F}_0:=\varphi(s)F(s,\theta_0)$. 
    Assume that Hypotheses~\ref{thm:root:ii}-\ref{thm:root:iv} 
    of Theorem~\ref{thm:1} holds, so that in particular, 
    $F(s,\cdot)$ is invertible around $\theta_0$.
    Define
    \begin{equation*}
	\cE_{\theta_0}F(s,x)
	:=
	F^{-1}\Paren*{s,F(s,\theta_0)\Paren*{1-\frac{x}{\varphi(s)}}}-\theta_0.
    \end{equation*}
    Let $\alpha_0:=0$ and $\alpha_k$ be the $k$-th Taylor coefficient of 
    $\cE_{\theta_0}F(s,\cdot)$ for $k\geq 1$,
    that is
    \begin{equation*}
	\alpha_k(s)
	:=\frac{\dD^k \cE_{\theta_0}F(s,0)}{k!}
	=
	\frac{\dD^k F^{-1}(s,F(s,\theta_0))}{k!}\left(-\frac{\coef{F}_0(s)}{\coef{F}_1(s)} \dD F(s,\theta_0) \right)^{\! k}.
    \end{equation*}
    Then
	$\Set{\alpha_k(s)}_{k=0}^p$ is $(\bbeta^{\upflat},p)$-related to $\bdelta^F(s)$, %$\Set{-\frac1{k!} \frac{\coef{F}_k(s)}{\coef{F}_1(s)}}_{k=0}^p$, 
	it satisfies equation~\eqref{eq:rs:s} (i.e.~Item~\ref{thm:root:i} of Theorem~\ref{thm:1}), and
    \begin{multline*}
	\theta(s)
	%=F^{-1}(s,0)
	= 
	\theta_0+\cE_{\theta_0}F(s,\varphi(s))
	\\
	=\theta_0+\cT_p \cE_{\theta_0}F(s,\varphi(s))+\grandO(\varphi(s)^{p+1})
	=
	\theta_0+\sum_{k=1}^p \alpha_k(s)\varphi(s)^k + \grandO(\varphi(s)^{p+1}).
    \end{multline*}
\end{proposition}
%%

%% PROOF
\begin{proof}
    By construction, for any $x$ small enough, 
    \begin{equation} \label{prop:inversion:eq}
	\frac{F(s,\theta_0+\cE_{\theta_0}F(s,x))-F(s,\theta_0)}{\dD F(s,\theta_0)}
	=\frac{-F(s,\theta_0)}{\dD F(s,\theta_0)}\frac{x}{\varphi(s)}. 
    \end{equation}
    In particular, $\cE_{\theta_0}F(s,0)=0$ and 
    $F(s,\theta_0+\cE_{\theta_0}F(s,\varphi(s)))=0$
    so that the root $\theta(s)$ to $F(s,\theta)=0$
    is $\theta(s)=\theta_0+\cE_{\theta_0}F(s,\varphi(s))$. 
    
	We set $\Set{\delta_k}_k:= \bdelta$. 
	Then, inspired by Example~\ref{ex:2}, we set 
	\begin{equation*}
	\psi(s,z)=\frac{F(s,\theta_0)-F(s,\theta_0+z)}{\dD F(s,\theta_0)}.
	\end{equation*}
	Hence 
	$\cT_p \psi(s,z)=\sum_{k=1}^p \delta_k(s) z^k$ which is analytic.
    From~\eqref{prop:inversion:eq},
	\begin{equation*}
    \psi(s,\cE_{\theta_0}F(s,x))=-\delta_0(s)x.
     \end{equation*}
	Note that 
	\begin{equation*}
    \fp_p \left( \cT_p \psi(s, \cT_p \cE_{\theta_0}F(s,x)) \right)=  \psi(s,\cE_{\theta_0}F(s,x)) =-\delta_0(s)x.
     \end{equation*}
	This ensures by Eq.~\eqref{eq:upflat:2}, together with Lemmas~\ref{lem:taylor-compos}-\ref{lem:taylor-compos-analytic}, that the coefficients $\balpha(s)$ of the Taylor expansion of $\cE_{\theta_0}F(s,\cdot)$
    around $0$ are 
    $(\bbeta^{\upflat},p)$-related to $\bdelta^F(s)$.
    Hypothesis~\ref{thm:root:i} is satisfied, so that 
    the conclusion of Theorem~\ref{thm:1} holds.
\end{proof}
%%

%%%%%%%%%%%%%%%%%%%%%%%%%%%%%%%%%%%%%%%%%%%%%%%%%%%%%%%%%%%%%%%%%%%%%%
\subsection{Perturbation of the roots and \texorpdfstring{$(\bbeta,p)$}{(β,p)}-related sequences} \label{sec:seq:pert}

Let $\bbeta \in \Set{\bbeta^{\upflat}, \bbeta^{\zig}}$. 
Let $\Lambda:\bS\to\RR$ and set $\coef{\Lambda}_0(s):=\varphi(s)^{\beta_0}\Lambda(s)$.
Let us consider solving 
\begin{equation*}
    F(s,\theta_\Lambda(s))=\Lambda(s).
\end{equation*}
We define 
\begin{equation*}
\coef{F_\Lambda}_k(s)
:=
\begin{cases}
    \coef{F}_0(s)-\coef{\Lambda}_0(s)&\text{ for }k=0,\\
    \coef{F}_k(s)&\text{ for }k>0
\end{cases}
\end{equation*}
as well as the corresponding $\delta_{\Lambda\,k}(s)$  and $\alpha_{\Lambda\,k}(s)$. 
%%

%% PROPOSITION
\begin{proposition}
    \label{prop:perturbation}
    For each $k\geq 1$, the coefficients
    $\alpha_{\Lambda\,k}(s)$ are polynomial expressions in 
    $\coef{F}_\ell(s)/\coef{F}_1(s)$, $\ell\leq k$ and $\coef{\Lambda}_0(s)/\coef{F}_1(s)$. 
    Moreover, if $\coef{\Lambda}_0/\coef{F}_1$ converges to $0$ with~$s$
    and provided that $\alpha_k$ converges, 
    then $\alpha_{\Lambda\,k}$ and $\alpha_k$ have
    the same limit for each $k\geq 0$.
\end{proposition}
\begin{proof}
The proof is immediate once it has been noticed by
$\delta_{\Lambda\,k}=\delta_k$ for $k>1$
and $\delta_{\Lambda\,0}=\delta_0-\coef{\Lambda}_0$.
The expression of $\alpha_{\Lambda\,k}(s)$ in terms of 
$\coef{F}_\ell(s)/\coef{F}_1(s)$ and $\Lambda(s)/\coef{F}_1(s)$
are obtained by the Newton formula. 
\end{proof}
%%

%%%%%%%%%%%%%%%%%%%%%%%%%%%%%%%%%%%%%%%%%%%%%%%%%%%%%%%%%%%%%%%%%%%%%%
%%%%%%%%%%%%%%%%%%%%%%%%%%%%%%%%%%%%%%%%%%%%%%%%%%%%%%%%%%%%%%%%%%%%%%
%%%%%%%%%%%%%%%%%%%%%%%%%%%%%%%%%%%%%%%%%%%%%%%%%%%%%%%%%%%%%%%%%%%%%%
%%%%%%%%%%%%%%%%%%%%%%%%%%%%%%%%%%%%%%%%%%%%%%%%%%%%%%%%%%%%%%%%%%%%%%

\section{Reparametrization and Changes of scale }

\label{sec:change-var-scale}

\subsection{Reparametrization through changes of variables}

\label{sec:change-var}

In the context of statistical estimation, it is sometimes 
convenient to perform some change of variable on the space
of parameters. 

Let us consider a one-to-one map $\Phi$ from the space of parameters $\Theta$
to the space~$\Theta^\dag$ with inverse $\Psi$.
We assume that $\Phi\in\cC^{p+1}$.
For the sake of simplicity,
we assume that~$0$ belongs to the interiors of both $\Theta$ and $\Theta^\dag$. 
Moreover, we assume that~$\Phi(0)=0$.

\begin{notation}[Pullback]
    For a function $f:\Theta\to\RR$, we set $\Psi^\star f:=f\circ \Psi$
    so that~$\Psi^\star f$ is a function from $\Theta^\dag$ to $\RR$.
    The operator $\Psi^\star$ is the \emph{pullback}. For $f:\bS\times\Theta\to\RR$, 
    we extend the pullback to $\Psi^\star f(s,\eta)=f(s,\Psi(\eta))$
    for $(s,\theta)\in\bS\times\Theta^\dag$.
	For a function $F\in\cC^{p+1}(\Theta,\RR)^\bS$, we define
	$F^\dag(s,\eta):=\Psi^\star F(s,\eta)=F(s,\Psi(\eta))$ which is a function 
	in $\cC^{p+1}(\Theta^\dag,\RR)^\bS$. 
	We also set
	$\coef{F^\dag}_k(s):=\varphi(s)^{\beta_k}\dD^k F^\dag(s,0)$
	for $k=0,\dotsc,p$.
\end{notation}

We focus on the case \textquote{up-flat} ($\bbeta=\bbeta^{\upflat}$).
From the Faà di Bruno formula (Proposition~\ref{prop:fdb}), 
we see that $\Set{\coef{F^\dag}_k(s)}_{k=0,\dotsc,p}$ 
converges whenever $\Set{\Set{F}_k(s)}_{k=0,\dotsc,p}$ does. 
This is not necessarily
true in the case \textquote{zig-zag}.
Indeed, let $\varphi(s)\in (-1,1)$, then in the \textquote{zig-zag} case, $\varphi(s)^{-1} \Set{F}_1(s) $ may diverge and so 
\[\Set{F^\dag}_2(s) = \Set{F}_2(s) (\dD \Psi(0))^2 + \varphi(s)^{-1} \Set{F}_1(s) \dD^2 \Psi(0)\] 
does not converge neither.
The same happens for $\Set{F^\dag}_{2k}(s)$, by Faà di Bruno formula~\eqref{eq:fdb} (in Proposition~\ref{prop:fdb}). So, starting from \textquote{zig-zag} case, taking $F^{\dag}$ we bump into the \textquote{up-flat} case.
More generally, if $\varphi(s)\in (-1,1)$, starting from $\bbeta=(\beta_0,\beta_1,\beta_2,\ldots,\beta_p)$, 
we bump into $(\beta_0,\beta_1,\beta_2^{\star},\beta_3^{\star},\ldots, \beta_p^{\star})$ 
where $\beta_k^\star:=\max_{n\in \{1,\ldots,k\}} \beta_n$. 
This latter kind of coefficient, as $\bbeta^{\upflat}$, is \textquote{stable} under such changes of variables.

\begin{proposition} \label{prop:change:expansion}
    Fix $s\in\bS$.
    Let us define $\bdelta=\Set{\delta_k}_{k=0,\dotsc,p}$ and $\bdelta^\dag=\Set{\delta_k}_{k=0,\dotsc,p}$ by
    \begin{equation*}
	\delta_k(s):=\frac{\coef{F}_k(s)}{k!}
	\text{ and }\delta^\dag_k(s):=\frac{\coef{F^\dag}_k(s)}{k!}
	\text{ for }k=0,1,\dotsc,p. 
    \end{equation*}
    Let $\balpha(s)$ (resp.~$\balpha^\dag(s)$) be the unique
    $(\bbeta^{\upflat},p)$-related sequence to $\bdelta(s)$ (resp.~$\bdelta^\dag(s)$)
    given in Lemma~\ref{lem:upflat}. 
    Then 
    \begin{equation} \label{eq:prop:change}
	 \alpha_q(s)
	=
	\sum_{m=1}^q \frac{\dD^m\Psi(0)}{m!}
	\sum_{k_1+\dotsb+k_m=q} \alpha_{k_1}^\dag(s)\dotsb  \alpha_{k_m}^\dag(s)
	\text{ for }q=1,\dotsc,p, 
    \end{equation}
    or, in other words, 
    \begin{equation*}
	\sum_{k=1}^p \alpha_k(s)^\bullet z^k:=\fp_p\circ {\cT_p^{\dummy} \Psi} 
	\Paren*{
\sum_{k=1}^p \alpha^\dag_k(s)^\bullet z^k
	},
    \end{equation*}
where ${\cT_p^{\dummy} \Psi}
\in \RR\Paren*{\bdelta^{\dummy}(s)}\dbra{z}$ is the formal Taylor expansion of $\Psi$ (see Definition~\ref{def:formal:taylor}).
\end{proposition}

The meaning of this proposition is that expanding the root $\theta^\dag(s)$ 
of $F^\dag(s,\theta^\dag(s))=0$ in terms of $\varphi(s)$ using our procedure 
of Lemma~\ref{lem:upflat}
and then extending $\Psi(\theta^\dag(s))$ using a Taylor series
leads to the same expansion as the one given computing 
the root~$\theta(s)$ of $F(s,\theta(s))=0$ using Lemma~\ref{lem:upflat}.
This is summarized in the diagram of Figure~\ref{cd:variable}.

\begin{figure}
\begin{center}
    \begin{tikzcd}[row sep=large, column sep=large]
	{(F,\theta)} 
	\arrow[r,"\text{expansion}"] 
	\arrow[d,"\Psi^\star\times\Phi"]
	&
    {\sum_{k=1}^p \alpha^\bullet_kz^k}
    \arrow[r,"\text{\vphantom{p}evaluation}"] 
	& 
    {\theta_0+\sum_{k=1}^p \alpha_k\varphi^k}
    % 2nd line
    \\
	{(F^\dag,\theta^\dag)} 
	\arrow[r,"\text{expansion}"] 
	&
	{\sum_{k=1}^p \alpha^{\dag\bullet}_kz^k}
	\arrow[r,"\text{\vphantom{p}evaluation}"] 
    \arrow[u,"{\fp_{p}\circ {\cT_p^{\dummysmall} \Psi}
%\Psi_{{\dummy}}
}"]
	&
	{\theta_0^\dag+\sum_{k=1}^p \alpha^{\dag}_k\varphi^k}
\end{tikzcd}
\caption{\label{cd:variable} Change of variable.}
\end{center}
\end{figure}

\begin{proof}
	In this proof let $F_{\dummy}:={\cT_p^{\dummy} F}$ and $F^\dag_{\dummy}:={\cT_p^{\dummy} F^\dag}$ from Definition~\ref{def:formal:taylor}.
    By \eqref{eq:formal:6} and  Lemma~\ref{lem:upflat},
    the sequence $\balpha(s)^\bullet\in\RR(\bdelta^{\dummy})^p$ corresponding 
    to $\alpha(s)$ is the unique solution to
\begin{equation}
    \label{eq:formal:4}
    \fs_{\bbeta^{\upflat}}\circ \fp_p\circ 
    F_{\dummy}(s,P(z))=0
    \text{ with }P(z):=\sum_{k=1}^p \alpha_k(s)^\bullet z^k \text{ and } \alpha_k(s)^\bullet \in \RR\Paren{\bdelta(s)^{\dummy}}.
\end{equation}
Similarly, $\balpha^\dag(s)^\bullet$ is the unique solution to 
\begin{equation*}
\fs_{\bbeta^{\upflat}}\circ \fp_p\circ 
    F^\dag_{\dummy}(s,P^\dag(z))=0
    \text{ with }P^\dag(z):=\sum_{k=1}^p \alpha_k^\dag(s)^\bullet z^k=0
\end{equation*}
and $F^\dag_{\dummy}(z)=\sum_{k=0}^p \delta_k^\dag(s)^{\dummy} z^k \in \RR\bra{\bdelta^\dag(s)^{\dummy}}\dbra{z}$.

The $\alpha_k^\dag(s)^\bullet$'s 
are rational functions 
in $\RR\Paren*{\bdelta^\dag(s)^{\dummy}}$.
From the very definition, $F^\dag(s,\eta):=F(s,\Psi(\eta))$. 
By the Faà di Bruno formula (Proposition~\ref{prop:fdb}), $\coef{F^\dag}_k(s)$ is a
linear combination of the elements~$\coef{F}_\ell(s)$, $\ell=1,\dotsc,k$ and the derivatives
of~$\Psi$. Therefore, the $\alpha_k^\dag(s)^\bullet$'s and $\delta_k^\dag(s)^{\dummy}$'s can
be seen as rational functions in 
$\RR\Paren*{\bdelta(s)^{\dummy}}$. 
Thus $F^\dag_{\dummy}(z)$ and ${\cT_p^{\dummy} \Psi}\circ P^\dag(z)$ can be seen as elements of $\RR\bra{\bdelta(s)^{\dummy}}\dbra{z}$.

By Lemma~\ref{lem:taylor-compos},
\begin{equation*}
\fp_p \circ F^\dag_{\dummy} = F^\dag_{\dummy} = {\cT_p^{\dummy} (\Psi^\star F)}
= \fp_p \circ  ({\cT_p^{\dummy} \Psi})^\star F_{\dummy}.
\end{equation*}
Since ${\cT_p^{\dummy} \Psi}\circ P^\dag(z) \in \RR\bra{\bdelta(s)^{\dummy}}\dbra{z}$ and 
\begin{equation*}
   0=\fs_{\bbeta^{\upflat}}\circ \fp_p\circ F^\dag_{\dummy}(s,P^\dag(z))
    =
    \fs_{\bbeta^{\upflat}}\circ \fp_p\circ F_{\dummy}(s,{\cT_p^{\dummy} \Psi}\circ P^\dag(z)).
\end{equation*}
Hence ${\cT_p^{\dummy} \Psi}\circ P^\dag(z) $ solves~\eqref{eq:formal:4}, 
which is known to have a unique solution, $P(z)$. Thus $P(z)={\cT_p^{\dummy} \Psi}\circ P^\dag(z)$. 
The result follows by an application of the evaluation map.
\end{proof}
%%

%%%%%%%%%%%%%%%%%%%%%%%%%%%%%%%%%%%%%%%%%%%%%%%%%%%%%%%%%%%%%%%%%%%%%%

\subsection{Change of scale}

\label{sec:change-scale}

We now introduce another change of variable, called \emph{change of scale},
leading to a different expansion. The impact of this change of 
scale technique is shown in Section~\ref{sec:exp-fam:natural}
to obtain a \textquote{natural} expansion in the context 
of exponential families.

Let us take a one-to-one function $\Phi\in\cC^{p+1}(\Theta,\Theta^\dag)$
whose derivative $\Phi'$ never vanishes.
Such a function $\Phi$ is seen as a \emph{scale}.

\begin{notation}
    \label{not:deriv}
    For a differentiable function $f:\RR\to\RR$, we write
    $\dD_\Phi f=\dD f/\dD\Phi$, where $\dD$ is the derivative operator. 
\end{notation}

Since we solve $F(s,\theta)=0$, a natural operation consists in setting
\begin{equation*}
    \widecheck{F}(s,\theta):=
    \Phi^\sharp F(s,\theta):=\frac{F(s,\theta)}{\dD\Phi(\theta)}
\end{equation*}
and to solve $\widecheck{F}(s,\theta)=0$. The root $\theta(s)$
is the same for $F$ and $\widecheck{F}$. However, as discussed below, 
the expansion could be different. 

In the case of the MLE, as $F(s,\cdot)$ takes scalar values, 
there exists a function $G\in\cC^{p+2}(\Theta,\RR)^\bS$ such that $\dD G(s,\cdot)=F(s,\cdot)$.
Considering $\widecheck{F}(s,\theta)$ corresponds to a change of variable at the level of the function $G$. 
Indeed, with Notation~\ref{not:deriv}, 
\begin{equation*}
    \dD_\Phi G(s,\theta)=\frac{\dD G(s,\theta)}{\dD\Phi(\theta)}=\widecheck{F}(s,\theta)
    \text{ in }
    \cC^{p+1}(\Theta,\RR)^\bS.
\end{equation*}
We set $G^\dag(s,x):=G(s,\Psi(x))$ and $F^\dag(s,x):=F(s,\Psi(x))$
so that from standard computations, 
\begin{equation*}
\dD (G^{\dag}) = (\dD_\Phi G)^{\dag}
\text{ and }
F^\dag = \dD_\Psi (G^\dag).
\end{equation*}
Therefore, 
\begin{equation*}
    \widecheck{F}^\dag:= (\widecheck{F})^\dag= \dD(G^\dag)
    \text{ so that }
    \widecheck{F}= \Phi^\star\dD(\Psi^\star G).
\end{equation*}

We denote by 
$\widecheck{\balpha}(s)$ and $\widecheck{\balpha}^\dag(s)$ the unique 
$(\bbeta^{\upflat},p)$-related sequences to 
$\Set{\coef{\widecheck{F}}_k(s)/k!}_{k=0,\dotsc,p}$
and $\Set{\coef{\widecheck{F}^\dag}_k(s)/k!}_{k=0,\dotsc,p}$ defined by Lemma~\ref{lem:upflat}.

Proposition~\ref{prop:change:expansion} may be applied to compute 
$\widecheck{\balpha}^\dag(s)$ and then $\widecheck{\balpha}(s)$, 
instead of computing the derivative of $\widecheck{F}$.
We summarize the various relationships in Figure~\ref{cd:scale}.

\begin{figure}
\begin{center}
    \begin{tikzcd}[row sep=large, column sep=large]
	% line 1
	G^\dag
	\arrow[r,"\dD_\Psi"]
	& 
	{(F^\dag,\theta^\dag)}
	\arrow[r,"\text{expansion}"] 
	& 
	{\sum_{k=1}^p \alpha^{\dag\bullet}_kz^k}
	\arrow[d,"{\fp_{p}\circ {\cT_p^{\dummysmall} \Psi}}"]
	\arrow[r,"\text{\vphantom{p}evaluation}"] 
	&
	{\theta^\dag+\sum_{k=1}^p \alpha^\dag_k\varphi^k}
	\\
	% line 2 
	G
	\arrow[r,"\dD"]
	\arrow[d,"\mathrm{Id}"]
	\arrow[u,"\Psi^\star"]
	&
	{(F,\theta)} 
	\arrow[r,"\text{expansion}"] 
	\arrow[d,"\Phi^\sharp\times\mathrm{Id}"]
	\arrow[u,"\Psi^\star\times\Phi"]
	& 
    {\sum_{k=1}^p \alpha^\bullet_kz^k}
	\arrow[r,"\text{\vphantom{p}evaluation}"] 
	&
	{\theta+\sum_{k=1}^p \alpha_k\varphi^k}
    \\
    % line 3
	G
	\arrow[r,"\dD_\Phi"]
	\arrow[d,"\Psi^\star"]
	&
    {(\widecheck{F},\theta)} 
	\arrow[r,"\text{expansion}"] 
	\arrow[d,"\Psi^\star\times\Phi"]
	&
	{\sum_{k=1}^p \widecheck{\alpha}^\bullet_kz^k}
	\arrow[r,"\text{\vphantom{p}evaluation}"] 
	&
	{\theta+\sum_{k=1}^p \widecheck{\alpha}\varphi^k}
    \\
    % line 4
	G^\dag
	\arrow[r,"\dD"]
	&
    {(\widecheck{F}^\dag,\theta^\dag)} 
	\arrow[r,"\text{expansion}"] 
	&
	{\sum_{k=1}^p \widecheck{\alpha}^{\dag\bullet}_k \varphi^k}
    \arrow[u,"{\fp_{p}\circ {\cT_p^{\dummysmall} \Psi}}"]
	\arrow[r,"\text{\vphantom{p}evaluation}"] 
	&
	{\theta^\dag+\sum_{k=1}^p \widecheck{\alpha}^\dag_k\varphi^k}
\end{tikzcd}
\caption{\label{cd:scale} Change of scales and change of variable.}
\end{center}
\end{figure}

Are $P_{\widecheck{\balpha}}$ and $P_{\balpha}$ the same?

Let us reduce to the case $\Phi(0)=0$. If $\Phi(x)=\lambda x$ for some $\lambda
\in \RR\setminus\Set{0}$, then the two expansions coincide by the fact that $F=
\dD G = \lambda \dD_\Phi G$ and Remark~\ref{rem:scalar:multipl}.

Note that if
$\alpha_1\neq 0$ and 
$\widecheck \alpha_1 \neq 0$ then they are different unless $\Phi''(0) = 0$ since
\begin{equation*}
    \widecheck \alpha_1(s) 
    =
    \frac{\alpha_1(s)}{1+\dfrac{\Phi''(0)}{\Phi'(0)}\alpha_1(s) \varphi(s)^{\beta_1-\beta_0}}.
\end{equation*}
Let us observe that $\widecheck{\alpha}_1$ and $\alpha_1$ coincide asymptotically as $s\to\infty$, when $\varphi(s)$ vanish and $\alpha_1(s)$ converges in distribution.

%%%%%%%%%%%%%%%%%%%%%%%%%%%%%%%%%%%%%%%%%%%%%%%%%%%%%%%%%%%%%%%%%%%%%%
%%%%%%%%%%%%%%%%%%%%%%%%%%%%%%%%%%%%%%%%%%%%%%%%%%%%%%%%%%%%%%%%%%%%%%
%%%%%%%%%%%%%%%%%%%%%%%%%%%%%%%%%%%%%%%%%%%%%%%%%%%%%%%%%%%%%%%%%%%%%%
%%%%%%%%%%%%%%%%%%%%%%%%%%%%%%%%%%%%%%%%%%%%%%%%%%%%%%%%%%%%%%%%%%%%%%
\section{Illustrating the results with some parametric models}

\label{sec:application}

We illustrate how to apply our results in this section focusing on the standard case of independent and identically distributed samples distributed according to a law of an exponential family 
and on an example for parameter estimation of stochastic processes.
Our results have been applied to provide new results on parameter estimation of Skew Brownian motion in~\cite{lm22}, 
where we also prove a phase transition as the one considered in this section for the binomial model (see Section~\ref{sec:symmetric}).

\subsection{Exponential family}

\label{sec:exp-fam}

In this section, we consider the class univariate models from the exponential
family of laws. This is a broad class of random variables~\cite{brown} with a
wide range of applications in statistical theory.

We consider an univariate parametric model from an exponential family of laws, 
which covers a wide variety of models~\cite{MR3221776}.
Let $\Theta$ be an open interval, and $T:\RR\to\RR$, $h:\RR\to\RR$ and $w,A \in \cC^{p+1}(\Theta)$, $p\geq 1$
be suitable functions. 
We consider $n$ independent copies $X_1,\dotsc,X_n$ with joint distribution under 
$\PP_\theta$ of the form 
\begin{equation}
    \label{eq:exp-fam:1}
	p(x_1,\ldots, x_n,\theta) = \left(\prod_{j=1}^n h(x_j) \right)  
	\exp\Paren*{w(\theta) \sum_{j=1}^n T(x_j) - n A(\theta)}.
\end{equation}
We consider that the functions $T$, $A$, $h$ and $w$ are such that the MLE exists
and the model is regular in particular
$\partial_\theta^j p(x,\theta)$ is integrable and the
following integrals are finite and equal $\int \partial_\theta^j p(x,\theta)
\vd x = \partial_\theta^j \int p(x,\theta) \vd x $, for $j=1,2$.

The MLE $\mlen$ is the zero of the score function 
\begin{equation*}
	\theta \mapsto F(n,\theta):= n \left( w'(\theta) T_n -  A'(\theta) \right) 
	\quad \text{where} \quad  
	T_n:=\frac1n \sum_{j=1}^n T(X_j).
\end{equation*}
Hence, $\mlen$ solves 
\begin{equation*} 
    \dD_w A(\mlen) := \frac{ A'(\mlen) }{ w'(\mlen) } = T_n,
\end{equation*} 
where Notation~\ref{not:deriv} is used. 

As the model is regular, $\EE_\theta \Prb{F(n,\theta)}= 0$. From the very definition of $F$, 
\begin{equation*}
\EE_\theta \Prb{F(n,\theta)} 
= n w'(\theta) \EE_\theta\Prb{T_n} - n A'(\theta)
\text{ and then }\EE_\theta\Prb{T_n}=\dD_w A(\theta),\ \forall \theta\in\Theta. 
\end{equation*}
The variance of $T_n$ may be computed as well:
\begin{equation}
    \label{eq:exp-fam:4}
	\Var\Prb{T_n} =\frac1{n}  \frac{A''(\theta) w'(\theta)- A'(\theta) w''(\theta)}{(w'(\theta))^3} = \frac1n \dD_w^2 A(\theta).
\end{equation}
In particular, the MLE $\mlen$ solves 
\begin{equation*}
    \EE_{\mlen}\Prb{T_n}=T_n.
\end{equation*}

If $T(X_1)$ is a square integrable random variable then the strong LLN and the CLT ensure respectively that
$T_n$ converges $\PP_\theta$ a.s.~to $\EE_\theta\Prb{T(X_1)}$
and, under $\PP_\theta$,
\begin{equation}
    \label{eq:exp-fam:2}
G_n:= \frac{T_n - \EE_\theta\Prb{T_n}}{\sqrt{\Var\Prb{T_n}}} 
    = \sqrt{n} \frac{T_n - \dD_w A(\theta)}{\sqrt{\dD_w^2 A(\theta)}}
    \xrightarrow[n\to\infty]{\mathrm{dist.}} G\sim\mathcal{N}(0,1).
\end{equation}

With \eqref{eq:exp-fam:4}, we define $\cI(\theta):=\dD_w^2 A(\theta) (w'(\theta))^2$.
This is the Fisher information of the model.

Taking $\bS=\NN$, $\varphi(s)=\frac1{\sqrt{s}}$, and $\bbeta=\bbeta^{\upflat}$, 
\begin{gather*}
    \coef{F}_0(n) 
    = \sqrt{n} \left(w'(\theta) T_n - A'(\theta)\right) 
    = \sqrt{\cI(\theta)} G_n 
     \cvdist[n\to\infty] \sqrt{\cI(\theta)} G =: \mu_0,
\end{gather*}
     and for $k=1,\dotsc,p$, 
\begin{gather*}
    \coef{F}_k(n) = w^{(k+1)}(\theta) T_n - A^{(k+1)}(\theta)  
    \cvprobat[n\to\infty] w^{(k+1)}(\theta) \dD_w A(\theta) - A^{(k+1)}(\theta) =: \mu_k.
\end{gather*}

By Proposition~\ref{prop:limit:alpha}, we get
the coefficients $\delta_k(n):=-\coef{F}_k(n)/k!\coef{F}_1(n)$, $k\geq 0$
of~$\bdelta(n)$ as well as 
the coefficients $\delta_k:=-\mu_k/k!\mu_1$ of its limit $\bdelta$. 
Theorem~\ref{th:random} (with $\theta_0=\theta$) provides the two approximations of the MLE estimator
$\mlen$ of $\theta$:
\begin{equation*}
	\theta_p(n) :=\theta + \sum_{k=1}^p \alpha_k(n) n^{-k/2} \qquad \text{and} \qquad
	\theta_p(n,\infty):= \theta + \sum_{k=1}^p \alpha_k n^{-k/2} 
\end{equation*}
with $\alpha_k(n)$ and $\alpha_k$ the $(\bbeta^{\upflat},p)$-related sequence to 
$\bdelta(n)$ and $\bdelta$ respectively (see~\eqref{eq:upflat:eval}).
Since $\mu_1=-\cI(\theta)$, 
\begin{equation*}
    \alpha_1(n) = - \frac{\coef{F}_0(n)}{\coef{F}_1(n) } 
	\cvdist[n\to\infty] - \frac{\mu_0}{\mu_1}
	= \frac{G}{\sqrt{\cI(\theta)}}.
\end{equation*}

Let us explicit the expansions with $p=3$:
\begin{equation}
	\label{eq:theta:expansion}
	\theta_3(s) = \theta + \frac{\alpha_1(n)}{\sqrt{n}}+  \delta_2(n) \left(\frac{\alpha_1(n)}{\sqrt{n}}\right)^{\!2}   + \left( 2  \delta_2(n) ^2 + \delta_3(n) \right) \left(\frac{\alpha_1(n)}{\sqrt{n}}\right)^{\!3},
\end{equation}
and
\begin{equation*}
    \theta_3(n,\infty) = \theta - \frac{G}{\sqrt{n\cI(\theta)}}
    -\frac{\mu_2}{2\mu_1} \left(\frac{G}{\sqrt{n\cI(\theta)}}\right)^{\!2}   
    + \left( \frac{1}{2}  \frac{\mu_2^2}{\mu^2_1}  -\frac{\mu_3}{6\mu_1} \right)
    \left(\frac{G}{\sqrt{n\cI(\theta)}}\right)^{\!3}.
\end{equation*}

Let us provide some examples of exponential families of laws to which the above procedure applies.

%% EXAMPLE
\begin{example}[Exponential distribution] 
    \label{exa:exponential}
We consider the Exponential distribution. The parameters is $\theta > 0$.
The related functions are 
\begin{gather*}
    T(x)=x,\ 
    w(\theta)=-\theta,\ 
    A(\theta)=-\log(\theta),
    \\ \text{so that }
    \dD_w A(\theta)=\frac{1}{\theta},\
    \dD^2_w A(\theta)=\cI(\theta)=\frac{1}{\theta^2}, \
    \coef{F}_{k}=\frac{(-1)^{k} k!}{ \theta^{k+1}} , k\geq 1.
\end{gather*} 
In this case the MLE is $\mlen=1/T_n$. The two approximations $\theta_p(n)$ and $\theta_p(n,\infty)$ in Theorem~\ref{th:random} are (using computations similar to the ones of Example~\ref{ex:1}) 
\begin{equation*}
    \theta_p(n)=\theta \sum_{k=0}^p \frac{(-G_n)^k}{n^{k/2}} 
\quad \text{ and } \quad
\theta_p(n,\infty)=\theta \sum_{k=0}^p \frac{(-G)^k}{n^{k/2}}
\end{equation*}
with $G_n=\sqrt{n}(\theta T_n-1)$ (See~\eqref{eq:exp-fam:2}) and $G\sim\cN(0,1)$.
With $p=+\infty$ since all the involved functions are analytic\footnote{The result may be found directly by using the relation between $\theta(n)$, $T_n$ and $G_n$.}, 
\begin{equation*}
    \theta_{\infty}(n)=\frac{\theta}{1+\dfrac{G_n}{\sqrt{n}}}
    =\theta\cH\Paren*{\dfrac{G_n}{\sqrt{n}}}
    =\theta + \theta H\Paren*{\dfrac{G_n}{\sqrt{n}}}
\end{equation*}
with 
\begin{equation*}
    \cH(x):=\frac{1}{1+x}\text{ and }H(x)=\frac{-x}{1+x}. 
\end{equation*}
Remark that $G_n>-\sqrt{n}$ almost surely, 
while $\theta_{\infty}(n,\infty)=\cH(G/\sqrt{n})$ on the event $\Set{G>-\sqrt{n}}$.
We then get the approximation of the distribution function of 
\begin{equation*}
\varepsilon(n):=\sqrt{n}(\theta_\infty(n)-\theta)
=
\sqrt{n}\Paren*{\frac{1}{T_n}-\theta}
=
\sqrt{n}H\Paren*{\frac{G_n}{\sqrt{n}}}
\approx 
\sqrt{n}H\Paren*{\frac{G}{\sqrt{n}}}
\end{equation*}
as 
\begin{equation*}
    \PP\Prb{\varepsilon(n)\leq x}
    \approx 
    \PP\Prb{\sqrt{n}(\theta_\infty(n,\infty)-\theta)\leq x}
    =
    \PP\Prb*{G\leq \frac{-x}{\theta+\frac{x}{\sqrt{n}}}},
\end{equation*}
using the Gaussian approximation of $G_n$ by $G$. 

With this example, we illustrate Fact~\ref{fact:5} in Section~\ref{sec:IFT}. 
We write $\alpha_1=\theta G$
with $G\sim\cN(0,1)$, as the Fisher information $\cI(\theta)$
is $1/\theta^2$.
In Figure~\ref{fig:exponential}, we represent
for several values of the sample size $n$ 
the following Kolmogorov-Smirnov
distances\footnote{We actually neglect the events $\Set{G\leq-\sqrt{n}}$ which 
is exponentially small.}, 
 \begin{equation}
     \label{eq:KSdist:1}
     \begin{split}
	 \Delta_1&:=\Delta(\varepsilon(n),\theta G),
 \
 \Delta_2:=\Delta(\epsilon(n),\theta\sqrt{n}H(G/\sqrt{n}))
 \\
 \text{ and }
	 \Delta_3&:=\Delta(\theta\sqrt{n}H(G/\sqrt{n}),\theta G)
 \end{split}
 \end{equation}
 with $\Delta$ defined in \eqref{eq:Kolm}. 
These distances are computed using~$10^6$ samples
of~$T_n$. The Kolmogorov-Smirnov statistics allows one to 
asserts that the Monte Carlo error, of order $10^{-3}$, is small in front
of the $\Delta_i$.
Since $H$ is invertible,  $\Delta_2=\Delta(G_n,G)$,
up to neglecting $\Set{G\leq -\sqrt{n}}$)
We then see that $\Delta_2$
is around half $\Delta_3$ (which does not depend on $\theta$).
This means that when $n$ is small, it is far
more satisfactory to use $\sqrt{n}H(G/\sqrt{n})$
than $G$ to construct for example confidence intervals or tests
for the value $\theta$.
\begin{figure}
\begin{center}
    \includegraphics{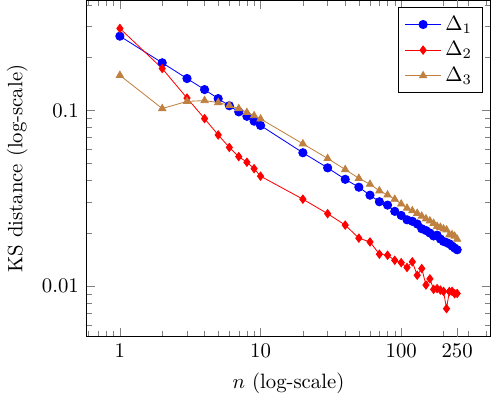}
    \caption{\label{fig:exponential} Kolmogorov-Smirnov distances defined
    in \eqref{eq:KSdist:1} for $\theta=1$.}
\end{center}
\end{figure}
\end{example}
%%

%% EXAMPLE
\begin{example}[Binomial distribution] \label{exa:binomial}
Let $\theta\in (0,1)$ be the parameter of a binomial model $(N,\theta)$ with $N\in \NN$ fixed ($N=1$ corresponds to the Bernoulli model).
Only the following functions differ from the previous case:
The related functions are 
\begin{gather*}
    T(x)=x,\ 
    w(\theta)=\log\frac{\theta}{1-\theta},\ 
    A(\theta)=-N \log(1-\theta)
    \\ 
    \text{so that }
    \dD_w A(\theta)=N \theta,\
    \dD^2_w A(\theta)=N \theta(1-\theta), \
	\cI(\theta)=\frac{N}{\theta(1-\theta)}.
\end{gather*}
In this case the MLE is the empirical mean divided by $N$. 
The derivative of the score and their a.s.~limit (LGN): 
for $m\geq 1$, 
\begin{equation*}
	\coef{F}_m(n) = n^{1-\beta_m/2} m!\left( (-1)^m \frac{T_n}{\theta^{m+1}} - \frac{(N-T_n)}{(1-\theta)^{m+1}}\right),
\end{equation*}
\begin{equation*}
	\coef{F}_m(\infty) := \mu_m = - m! N \left(  \frac{(-1)^{m+1}}{\theta^{m}} + \frac{1}{(1-\theta)^{m}} \right).
\end{equation*}
Besides, the CLT ensures that 
\begin{equation*}
   % \coef{F}_0(1/\sqrt{n})
	\coef{F}_0(n) = \sqrt{n} \frac{T_n- N \theta}{\theta(1-\theta)} = G_n \sqrt{\cI(\theta)} \cvdist[n\to\infty]
%    \coef{F}_0(0) 
	\coef{F}_0(\infty)  :=\sqrt{\cI(\theta)} G.
\end{equation*}
And $\coef{F}_1(\infty)= - \cI(\theta)$.

In particular the two approximations $\theta_p(n)$ and $\theta_p(n,\infty)$ (for $p=2$) in Theorem~\ref{th:random} are
\begin{multline*}
\theta_2(n)=\theta+\frac{G_n}{\sqrt{n}}\frac{\sqrt{N \theta^3 (1-\theta)^3}}{T_n(1-2\theta)+N\theta^2}
    +\frac{G_n^2}{n} N {\theta^2 (1-\theta)^2} \frac{\Paren{\theta^3+(1-\theta)^3} T_n - \theta^3 N}{\Paren{ (1-2\theta) T_n + N \theta^2 }^3 }.
\end{multline*}
and
\begin{equation*}
    \theta_2(n,\infty)=\theta + \frac{G\sqrt{\theta(1-\theta)}}{\sqrt{n\cdot N}}   
        +\frac{G^2}{n}\frac {1-2 \theta}{ N}.
\end{equation*}
\end{example}

%% BOUNDARY LAYER
To illustrate a \textit{boundary layer} effect, let us consider $\theta_3(n,\infty)$ in the previous example:
\begin{equation*}
    \theta_3(n,\infty)=\theta + \frac{G\sqrt{\theta(1-\theta)}}{\sqrt{n\cdot N}}   
        +\frac{G^2}{n} \frac {1-2 \theta}{ N}
			+ \frac{G^3}{(n\cdot N)^{3/2}} \frac{1 -5\theta +5 \theta^2}{\sqrt{\theta(1-\theta)}}.
\end{equation*}			
If $\theta$ is close to $0$ or to $1$ then the coefficient of $\frac{1}{\sqrt{n}}$ goes to $0$ and 
the other ones give their maximal contribution. In particular the coefficient of $\frac{1}{n^{3/2}}$ goes to infinity as $\theta$ goes to $0$ or to $1$.
Actually the expansion we used is not adapted when $\theta$ is close to~$0$ or~$1$.
This is  a \emph{boundary layer}.
More precisely, assume for simplicity $N=1$, $n=10^{m}$, and $\theta=10^{-t}$ (or $\theta=1-10^{-t}$) for $m,t\in \NN$:
$\theta_3(n,\infty) -\theta = 10^{-(m+t)/2} G + 10^{-m} G^2 + 10^{-(3m-t)/2} G^3$.\\
If $t=m$ all terms provide a contribution of the same order $10^{-m}$.
\\
If $t>m$, say $t=2m$ then
$\mlen -\theta \approx 10^{-3m/2} G + 10^{-m} G^2 + 10^{-m/2} G^3$ and the dominant term is the last one. 

%% REMARK
\begin{remark}
Note that if $\theta$ is close to $1/2$ then the coefficient of $\frac{1}{\sqrt{n}}$ in $\theta_3(n)$ gives its maximal contribution, the one of $\frac{1}{n}$ is close to $0$, and the coefficient of $\frac{1}{n^{3/2}}$ gives its maximal contribution (i.e.~$\frac{-1}{2\sqrt{N}} G^3$).
Moreover the latter term becomes $0$ for $\theta = \frac{5\pm \sqrt 5}{10} \in (0,1)$. % which are  lose to $0.27$ and $0.73$. 
The case $\theta = 1/2$ allows for a different choice of $\bbeta$: $\bbeta=\bbeta^{\zig}$. This is shown in Section~\ref{sec:symmetric}.
\end{remark}
%%

%%%%%%%%%%%%%%%%%%%%%%%%%%%%%%%%%%%%%%%%%%%%%%%%%%%%%%%%%%%%%%%%%%%%%%
\subsection{Using the natural parameters for the exponential family}

\label{sec:exp-fam:natural}

We saw in Sections~\ref{sec:change-var} and \ref{sec:change-scale}
two ways to consider another parametrization for the space of parameters. 
Such parametrizations could lead to simpler expressions. 

In the case of the exponential family, let us discuss the use
of the \textquote{natural scale}. For this, we
use the change of variable $\Phi:=w$ (with inverse $\Psi$), 
where the notations of~\eqref{eq:exp-fam:1} are in use. 

We consider that $A$ and $w$ are of class $\cC^{p+1}$ (or are analytic, 
in which case~$p=+\infty$). 

The likelihood and $\widecheck{F}$ are
\begin{equation*}
    \cL_n(\theta)=n(w(\theta)T_n-A(\theta))
   \text{ and }
    \widecheck{F}(n,\theta)=n(T_n-\dD_w A(\theta)).
\end{equation*}
The MLE $\mlen$ solves $\mlen = \dD_w A^{-1}(T_n)$ that can be rewritten as
\begin{equation*}
   \mlen =\cH_{\theta_0}\Paren*{\frac{G_n}{\sqrt{n}}}
    \text{ with }
    \cH_\theta(x):=\dD_w A^{-1}\Paren*{\dD_w A(\theta)+x\sqrt{\dD^2_w A(\theta)}}
\end{equation*}
with $G_n$ defined by \eqref{eq:exp-fam:2}.
Proposition~\ref{prop:inversion} shows that
\begin{equation*}
\theta_p(n)=\cT_p  \cH_{\theta_0} ({G_n}/{\sqrt{n}})
\end{equation*}
that is
\begin{equation}
    \label{eq:exp-fam:3}
\theta_p(n)
=\theta_0+\sum_{k=1}^p\widecheck{\alpha}_k(n)n^{-k/2}
=\theta_0+\sum_{k=1}^p \omega_k \Paren*{\frac{G_n}{\sqrt{n}}}^k 
\end{equation}
    with 
    \begin{equation*}
	\widecheck{\alpha}_k(n):=G^k_n \omega_k
	\text{ and }
	\omega_k :=
    \frac{\dD^k(\dD_w A^{-1})(\dD_w A(\theta_0))}{k!}
    \Paren{\dD^2_w A(\theta_0)}^{k/2}.
    \end{equation*}
    In particular, the $\omega_k$'s are deterministic and do not depend on $n$. 
    They are the Taylor coefficients of $\cH_{\theta_0}$.
%%

%% EXAMPLE
\begin{example}[Binomial distribution] \label{exa:binomial:2}
    For the binomial distribution (see Example~\ref{exa:binomial}), this expression simplifies for all $p\in \{1,\ldots,\infty\}$
    as 
    \begin{equation*}
	\theta_p(n)
	=\theta+\frac{\sqrt{\theta(1-\theta)}}{\sqrt{N\cdot n}}G_n
	= T_n
    \end{equation*}
	which is the empirical mean.
	Unlike for the exponential distribution,
    the distance between $\sqrt{n}(\theta_p(n)-\theta)$
    and the normal distribution $\cN(0,\theta(1-\theta)/N)$ comes
    solely
    from the distance between $G_n$ and the distribution $\cN(0,1)$.
\end{example}

Note that in Example~\ref{exa:exponential} the change of variable does not provide new approximations for the MLE because the model is basically already in natural scale. 

%%%%%%%%%%%%%%%%%%%%%%%%%%%%%%%%%%%%%%%%%%%%%%%%%%%%%%%%%%%%%%%%%%%%%%

\subsection{A phase transition of the approximation}

\label{sec:symmetric}

The possibility to change scale (\textit{i.e.},~$\bbeta=\bbeta^{\zig}$) in Theorem~\ref{th:random} sometimes leads to a phase transition in the approximation of the MLE estimator.
This happens when the distribution of the score is symmetric with respect to a parameter $\theta_0 \in \Theta$.

Let us consider the case of the binomial model in~Example~\ref{exa:binomial}.
If the true parameter~$\theta=\frac12$ then for all $m\in \NN$ even, when $\bbeta= \bbeta^{\upflat}$,
\begin{equation*}
    \coef{F}_{m}(n)= \frac{m!}{\sqrt{n} \theta^{2m}} \coef{F}_0(n) \cvprobat[n\to \infty] 0
\end{equation*}
If $\bbeta=\bbeta^{\zig}$ then the latter coefficients for $m$ even are
\begin{equation*}
\coef{F}_{m}(n)= \frac{m!}{\theta^{m}} \coef{F}_0(n) \cvdist[n\to \infty]  \mu_m:=\frac{m!}{\theta^{m}} \mu_0 = \frac{m!}{\theta^{m}} \sqrt{\cI(\theta)} G
\end{equation*}
and for $m$ odd, $\coef{F}_{m}(n)= - \frac{m! N}{\theta^{m+1}}= \frac{m!}{\theta^{m-1}} \coef{F}_{1}(n)$.

Applying Theorem~\ref{th:random} with $\bbeta=\bbeta^{\zig}$ and $\theta_0=\theta= \frac12$ provides new approximations of the MLE $\mlen$ (see Lemma~\ref{lem:zigzag} or equation~\eqref{eq:zig:eval}). 
The expansion for $n\in\NN$, $p\in \{1,\ldots, \infty\}$, coincides with the one obtained in the previous section:
\begin{equation*}
	\label{eq:theta:expansion:zigzag}
	\theta_p(n) = \theta + \frac{\alpha_1(n)}{\sqrt{n}}.
\end{equation*}
Indeed all coefficients $\alpha_k(n)=0$ for $k\geq 3$.  By Lemma~\ref{lem:zigzag}  (equivalently, equation~\eqref{eq:zig:eval}) it suffices to prove that $\alpha_{2k+1}=0$ for $k\geq 1$ and this can be shown by induction using the fact that $\coef{F}_{2 k}(n)= \frac{(2k)!}{\theta^{2k}} \coef{F}_0(n)$ and $\coef{F}_{2k+1}(n)= \frac{(2k+1)!}{\theta^{2k}} \coef{F}_{1}(n)$.

In the case of the MLE estimator of the skewness parameter of Skew Brownian motion, such a phase transition has been noticed in~\cite{lejay2014} and made explicit in~\cite{lm22}.

%%%%%%%%%%%%%%%%%%%%%%%%%%%%%%%%%%%%%%%%%%%%%%%%%%%%%%%%%%%%%%%%%%%%%%
%%%%%%%%%%%%%%%%%%%%%%%%%%%%%%%%%%%%%%%%%%%%%%%%%%%%%%%%%%%%%%%%%%%%%%
%%%%%%%%%%%%%%%%%%%%%%%%%%%%%%%%%%%%%%%%%%%%%%%%%%%%%%%%%%%%%%%%%%%%%%

\subsection{An application to statistics of diffusions}

\label{sec:diffusion}

We show how to use our method for non i.i.d. random 
variables and another estimator than the MLE.

Let us consider the inference of the parameter $\theta$
in the Stochastic Differential Equations (SDE)
\begin{equation*}
    \vd X^\theta_t = \theta\cdot b(X^\theta_t)\vd t+\sigma \vd B_t,\ t\geq 0. 
\end{equation*}
for a Brownian motion $B$, a function $b$ regular enough, and $\sigma \in (0,+\infty)$.

For this, we follow the methodology of \cite{Bibby1995}, which 
is based on the martingale approximations of the log-likelihood.

We fix a time step $\Delta>0$. We assume we observe a realization 
$(\sX_0,\dotsc,\sX_n)$ of 
the vector $(X_0^\theta,X_1^\theta,\dotsc,X_n^\theta)$
where $X^\theta_i$ is a shorthand for $X^\theta_{i\Delta}$, $i=0,\dotsc,n$.

The unknown parameter $\theta$ is estimated by 
$\mlen$ that solves $F_n(\mlen)=0$ with 
\begin{equation}
    \label{eq:sd:3}
    F_n(\theta):=\sum_{i=1}^n \frac{b(\sX_{i-1})}{\sigma^2}(\sX_i-M(\sX_{i-1},\theta))
\end{equation}
     where 
\begin{equation*}
    M(x,\theta):=\EE\Prb{X^\theta_{\Delta}\given X^\theta_0=x}. 
\end{equation*}

The asymptotic of $\mlen$ is studied when $n\to\infty$ while $\Delta$
is kept fixed (long time and fixed time step).

By a straightforward computation,
\begin{equation*}
    \dD_\theta^k F_n(\theta):=
    -
    \sum_{i=1}^n \frac{b(\sX_{i-1})}{\sigma^2}\dD_\theta^k M(\sX_{i-1},\theta).
\end{equation*}

%% HYPOTHESIS
\begin{hypothesis}
    \label{hyp:sd:1}
The process $X^{\theta}$ is ergodic. Its invariant measure 
is its renormalized speed measure
\begin{equation*}
    \mu_\theta(\vd x)
    =\frac{1}{A(\theta)}\exp\Paren*{2\theta\int_0^x \frac{b(y)}{\sigma^2}\vd y} \vd x,
\end{equation*}
where $A(\theta)$ is a normalizing constant such that $\int_\RR \mu_\theta(\vd x)=1$.
\end{hypothesis}

We denote by $(t,x,y)\mapsto p_\theta(t,x,y)$ the density transition function of $X^\theta$.
For $f:\RR^2\to\RR$, we write 
\begin{equation*}
    Q^\theta(f):=\int_\RR\int_\RR f(x,y)p_\theta(\Delta,x,y) \mu_\theta(\vd x).
\end{equation*}

The following proposition is a consequence of ergodicity.

%% PROPOSITION
\begin{proposition}[{See \textit{e.g.}, \cite[Lemma~3.1]{Bibby1995}}]
    Let $f:\RR^2\to\RR$ be a measurable function such 
    that $Q^\theta(f^2)<+\infty$. Then 
    \begin{equation*}
	\frac{1}{n}\sum_{i=1}^n f(X^\theta_{i-1},X^\theta_i)\cv[n\to\infty] Q^\theta(f) 
    \end{equation*}
    in $\mathrm{L}^2(\PP_\theta)$. If $Q^\theta(f)=0$, then 
    \begin{equation*}
	\frac{1}{\sqrt{n}}\sum_{i=1}^n f(X^\theta_{i-1},X^\theta_i)\cvdist[n\to\infty] 
	G\text{ with }G\sim\cN(0,Q^\theta(f^2))
    \end{equation*}
    under $\PP_\theta$.
\end{proposition}

Hence, for $k\geq 1$, owing to the definition of the invariant measure, 
under $\PP_\theta$, replacing the observations $\sX_i$ by the corresponding
random variable $X^\theta_i$, 
\begin{equation*}
    \frac{-\dD_\theta^k F_n(\theta)}{n}
    =
    \frac{1}{n}
    \sum_{i=1}^n \frac{b(X^\theta_{i-1})}{\sigma^2}\dD_\theta^k M(X^\theta_{i-1},\theta)
    \cvproba[n\to\infty]
    \int_{\RR} \frac{b(x)}{\sigma^2}\dD^k_\theta M(x,\theta)\mu_\theta(\!\vd x),
\end{equation*}
while, if  $\int_{\RR} \frac{b(x)}{\sigma^2}\dD^k_\theta M(x,\theta)\mu_\theta(\!\vd x)=0$,
\begin{equation*}
    \frac{1}{\sqrt{n}} \dD^k F_n(\theta) = \frac{1}{\sqrt{n}}
    \sum_{i=1}^n \frac{b(X^\theta_{i-1})}{\sigma^2}\dD_\theta^k M(X^\theta_{i-1},\theta)
    \cvdist[n\to\infty] \cN(0,c^2)
\end{equation*}
with 
\begin{equation*}
    c^2=\int_{\RR}\int_{\RR}
    p(\Delta,x,y)\frac{b(x)^2}{\sigma^4}(y-M(x,\theta))^2\mu_\theta(\dd x)\vd y.
\end{equation*}

We note that the quantity $M(x,\theta)$ cannot be always computed in an explicit manner, because
\begin{equation}
    \label{eq:sd:1}
    M(x,\theta)= \EE_\theta\Prb{X^\theta_\Delta \given X^\theta_0=x}=x+\theta\int_0^\Delta \EE_\theta\Prb{b(X_s^\theta)\given X^\theta_0=x}\vd s.
\end{equation}

%% EXAMPLE : ORNSTEIN-UHLENBECK
\begin{example}[Ornstein-Uhlenbeck process]
Let us consider that $b(x)=-x$ and $\theta>0$ so that $X^\theta$ is the Ornstein-Uhlenbeck process.
From \eqref{eq:sd:1}, 
\begin{equation*}
    M(x,\theta)=x\exp(-\Delta \theta)
    \text{ and }
    \mu_\theta(\dd x)=\frac{\sqrt{\theta}}{\sqrt{\pi \sigma^2}}\exp\Paren*{-\frac{\theta x^2}{\sigma^2}}\vd x
\end{equation*}
which is a centered Gaussian law of variance $\sigma^2/(2\theta)$.
Therefore, $\dD^k M(x,\theta)=x (-\Delta)^k\exp(-\Delta \theta)$. 
We are then in the case \textquote{up-flat}. 
We consider the change of variable $\eta:=\Phi(\theta):=e^{-\Delta \theta}$ (see Section~\ref{sec:change-var}) and we get
\begin{equation*}
	F^\dag_n(\eta) = 
	\sum_{i=1}^n
    \frac{-\sX_{i-1}}{\sigma^2}(\sX_i-\sX_{i-1}\eta).
\end{equation*}
It holds that 
$\dD_\eta F^\dag_n = \sum_{i=1}^n \frac{1}{\sigma^2}(\sX_{i-1})^2$ and $\dD_\eta^k F^\dag_n =0$.
By Theorem~\ref{th:random},
\begin{equation*}
    \eta(n)
	= \eta - \frac{F^\dag_n(\eta)}{\dD_\eta F^\dag_n(\eta)}
	= \frac{\sum_{i=1}^n \sX_{i-1}\sX_i}{\sum_{i=1}^n (\sX_{i-1})^2}
\end{equation*}
so that 
\begin{equation*}
    \mlen
=\frac{‐1}{\Delta}\log\frac{
    \sum_{i=1}^n \sX_{i-1}\sX_i
}{\sum_{i=1}^n (\sX_{i-1})^2}.
\end{equation*}
We have the same formula as the one in \cite[Example 2.1, p.~21]{Bibby1995}, 
which is obtained by direct means (note that we use $b(x)=-x$ instead of $b(x)=x$ as in \cite{Bibby1995}).
\end{example}
%%

%%% EXAMPLE: APPROXIMATION
%%
\begin{example}[Using an approximation of $M(x,\theta)$]
    When no close form solution to $M(x,\theta)$ (see \eqref{eq:sd:1}), 
    we consider that we have a smooth approximation $M_\delta(x,\theta)$
    of $M(x,\theta)$ depending on a parameter $\delta$. 
    We set 
    \begin{equation*}
	\epsilon_\delta(\theta):=\sum_{i=1}^n \frac{b(\sX_{i-1})}{\sigma^2}
	(M_\delta(\sX_{i-1},\theta)-M(\sX_{i-1},\theta))
    \end{equation*}
    and define $F_{\delta,n}$ similarly 
    to $F_n(\theta)$ in \eqref{eq:sd:3} with $M_\delta$ instead of $M$.
    The solution $\theta_{\delta,n}$ to $F_{\delta,n}(\theta_{\delta,n})=0$
    also solves
    \begin{equation}
	\label{eq:sd:4}
	F_n(\theta_{\delta,n})=\epsilon_{\delta}(\theta_{\delta,n}).
    \end{equation}
    Assume that we are in the case \textquote{up-flat} and 
    that $\epsilon_{\delta}(\theta_{\delta,n})/\sqrt{n}$ converges to $0$
    as $(\delta,n)$ converges to $(0,+\infty)$ as a net (this could imply 
    some relationship between the rate of convergence of $\delta$ to $0$
    and the one of $n$ to $+\infty)$.
    Proposition~\ref{prop:expansion} allows one to conclude
    that $\theta_{\delta,n}$ and $\mlen$ have the same 
    asymptotic behavior.
\end{example}
%%

%%%%%%%%%%%%%%%%%%%%%%%%%%%%%%%%%%%%%%%%%%%%%%%%%%%%%%%%%%%%%%%%%%%%%%

%%%%%%%%%%%%%%%%%%%%%%%%%%%%%%%%%%%%%%%%%%%%%%%%%%%%%%%%%%%%%%%%%%%%%%
%%%%%%%%%%%%%%%%%%%%%%%%%%%%%%%%%%%%%%%%%%%%%%%%%%%%%%%%%%%%%%%%%%%%%%
%%%%%%%%%%%%%%%%%%%%%%%%%%%%%%%%%%%%%%%%%%%%%%%%%%%%%%%%%%%%%%%%%%%%%%
%%% CONCLUSION 

\section{Conclusion}

We have given an expansion of the solution $\mlen$ to a problem 
$F_n(\mlen)=0$ in a way which involves the asymptotic behavior
of $F_n$ and its derivatives. This covers (quasi-, pseudo-) maximum likelihood
estimation and Generalized Method of Moments. We could also consider 
errors in the observations.
Remark that $n$ is a quantity related to the asymptotic behavior 
(\textit{e.g.},~the sample size). Our hypotheses are rather weak 
so that we may consider independent samples as well as dependent 
samples such as the ones arising in statistics of diffusion or chronological series.
With the inference in view, we have an expansion of type 
$\mlen=\theta_0+G_n\varphi(n)+E_n(G_n\varphi(n))+\grandO(\varphi(n)^{p+1})$
where $G_n$ converges in distribution to some random variable $G$, 
$\varphi(n)$ is a rate (say $\varphi(n)=1/\sqrt{n})$ and $E_n$ is 
a non-linear, random function with $E_n(x)=\grandO(x^2)$. 
In the cases considered, the convergence of $G_n$
depends on the one of $F_n$ and its first derivative. 
The function $E_n$ has a Taylor expansion whose coefficients 
have an asymptotic behavior that depends on the derivatives of $F_n$.

Our expansion can be used to access semi-asymptotic behaviors of the estimator, 
as we have illustrated in a toy example in Section~\ref{sec:exp-fam}.
The distribution of $\varphi(n)(\mlen-\theta_0)$ is close
to the one of $G$. Its departure from $G$ is due to the distance
between the distributions of $G_n$ and the one of $G$, as 
well as the behavior of $E_n$. We illustrated this in Fact~\ref{fact:5} (Section~\ref{sec:IFT}) 
and in Example~\ref{exa:exponential}.
So, Berry-Esséen type results seems accessible as well because (in the cases we specify) the extension we establish
is a non-linear transformation of an asymptotically pivotal quantity (see for instance~\cite{pinelis16a,pinelis17a}).
Note that Berry-Esséen results on the convergence of the latter to its limit may be combined with our expansion
(as in the case of i.i.d.~exponential random variable). Similar approaches should 
also apply to Wasserstein distance (see \textit{e.g.}, \cite{anastasiou17a,anastasiou17b}).

Let us mention that the main results of this document (in particular Theorem~\ref{thm:1}) 
can be easily extended also for multidimensional parameters 
thanks to multidimensional Inverse Function Theorem (see, \textit{e.g.},~\cite{abraham}).
It should also hold true in an infinite dimensional setting, 
which should be useful to deal with semiparametric models although 
there are some fine technical points to deal with (see the discussion in 
\cite{gill89a,gill89b,MR1652247}). In the multidimensional setting, 
a Taylor expansion is used in~\cite{zbMATH03551708} to compute the first two moments of functions of estimators. This might be combined easily with our approach to obtain higher 
moments. We expect to compute cumulants \cite{zbMATH07523688} as well in a wide variety of 
cases, and to identify some corrections to enforce the quality 
of the estimation~\cite{MR1652247}.

%%%%%%%%%%%%%%%%%%%%%%%%%%%%%%%%%%%%%%%%%%%%%%%%%%%%%%%%%%%%%%%%%%%%%%
%%%%%%%%%%%%%%%%%%%%%%%%%%%%%%%%%%%%%%%%%%%%%%%%%%%%%%%%%%%%%%%%%%%%%%
%%% APPENDIX %%%%%%%%%%%%%%%%%%%%%%%%%%%%%%%%%%%%%%%%%%%%%%%%%%%%%%%%%
%%%%%%%%%%%%%%%%%%%%%%%%%%%%%%%%%%%%%%%%%%%%%%%%%%%%%%%%%%%%%%%%%%%%%%
%%%%%%%%%%%%%%%%%%%%%%%%%%%%%%%%%%%%%%%%%%%%%%%%%%%%%%%%%%%%%%%%%%%%%%
\appendix

%%%%%%%%%%%%%%%%%%%%%%%%%%%%%%%%%%%%%%%%%%%%%%%%%%%%%%%%%%%%%%%%%%%%%%
\section{Asymptotic inversion}

\label{sec:asymopt-inversion}

The goal of the section is proving Theorem~\ref{thm:1}.
The proof of Theorem~\ref{thm:1:analytic} is similar, yet simpler 
as there is no need to control the remainder.

%%%%%%%%%%%%%%%%%%%%%%%%%%%%%%%%%%%%%%%%%%%%%%%%%%%%%%%%%%%%%%%%%%%%%%
\subsection{A Taylor formula with remainder}

In this section, we basically consider the $p$-th order Taylor polynomial of 
$z \mapsto F(s,\pP{p}(z;\balpha))$ around 0 and a convenient expression for the remainder.
\begin{hypothesis} \label{hyp:P}
    Let $p\in \NN\setminus\{0\}$.
    Let $\balpha:=\Set{\alpha_k}_{k=0,\dotsc,p}\in\RR^{p+1}$
    with $\alpha_0=0$.
With $\pP{m}$ defined by
\begin{equation}
    \label{eq:P}
    \pP{m}(z; \balpha) :=\sum_{k=0}^m \alpha_k z^k\text{ for }m=0,1,\dotsc,p
    \text{ and }z\in\RR, 
\end{equation}
we assume that $\pP{m}(\II;\balpha) \subseteq \Theta$ for all $m\in \Set{1,\ldots,p}$.
\end{hypothesis}
For the sake of simplicity, in this section, we consider a real-valued function $F$ on~$\Theta$,
but we could consider a real-valued function on $\bS \times \Theta$.
\begin{notation}
    \label{not:remainder}
	   For $m=1,\dotsc,p$, $z\in\RR$, and $F\in\cC^p(\Theta,\RR)$ we define 
\begin{multline}
\label{eq:def:Rm}
    R_m(z;\balpha):=
    \sum_{\mathclap{k_1+\dotsb+k_m\leq p}}
    \alpha_{k_1}\dotsb\alpha_{k_{m}} z^{k_1+\dotsb+k_{m}}
    \\
    \times 
    \int_0^1
    \int_0^{\tau_1}
    \dotsi
    \int_0^{\tau_{m-1}}
    \Big(
\dD^m F\Paren[\big]{\tau_m \pP{p- (k_1 + \dotsb +k_{m-1})} (z;\balpha) }
\\
-
\dD^m F\Paren[\big]{\tau_m \pP{p- (k_1 + \dotsb +k_{m})} (z;\balpha)  }
\Big)\vd\tau_m\dotsb\vd \tau_1.
\end{multline}
\end{notation}
\begin{notation}
    \label{not:J}
    We define for $z\in\II$, 
    \begin{align}
	\notag
	\pP{p}^\star(z;\balpha)&:=\sum_{k=1}^p \abs{\alpha_k}\cdot \abs{z}^k,
	\\
	\notag
	J(z;\balpha)&:=[-\pP{p}^\star(z;\balpha),\pP{p}^\star(z;\balpha)] \cap \Theta ,
    \\
	\label{eq:normF}
	\text{and }
	\abs{F}_m(z;\balpha)&:=\sup_{w\in J(z;\balpha)}\abs{\dD^m F(w)}.
    \end{align}
\end{notation}
\begin{lemma} \label{lem:devel:1}
Assume Hypothesis~\ref{hyp:P}.
Let $F\in\cC^{p+1}(\Theta,\RR)$. 
Then for all $z\in \II$, 
\begin{multline}
    \label{eq:102}
    F(\pP{p}(z;\balpha))-F(0)
    \\
    =\sum_{m=1}^p \frac{1}{m!}\dD^{m} F(0)  \sum_{\mathclap{k_1+\dotsb+k_m\leq p}}
    \alpha_{k_1}\dotsb\alpha_{k_{m}}
    z^{k_1+\dotsb +k_{m}}+\sum_{m=1}^p R_m(z;\balpha).
\end{multline}
Besides, for $m=1,\dotsc,p$, 
\begin{equation}
    \label{eq:ineq:Rm}
    \abs{R_m(z;\balpha)}
    \leq \abs{F}_{m+1}(z;\balpha)\pP{p}^\star(1;\balpha)^{m+1}\abs{z}^{p+1}
    \text{ for } \abs{z}\leq 1.
\end{equation}
\end{lemma}
%%%

%%%
\begin{proof}
For simplicity we do not specify the dependence on the sequence $\balpha$.
By the fundamental theorem of calculus,
\begin{multline}
    F(\pP{p}(z))-F(0)=\int_0^1 \dD F\Paren{\tau \pP{p}(z) }\cdot \pP{p}(z) \vd \tau
	\\
	=\dD F(0) \pP{p}(z)+ \int_0^1 ( \dD F\Paren{\tau_1 \pP{p}(z)}  - \dD F(0)) \cdot \pP{p}(z)\vd \tau_1.
\end{multline}
Adding and removing another suitable term yields
\begin{multline}
\label{eq:100}
F(\pP{p}(z))-F(0)  
	\\
	  =
	  \dD F( 0) \pP{p}(z)
    +
    \sum_{k_1=1}^p\alpha_{k_1} z^{k_1}
\int_0^1 \Paren[\Big]{\dD F\Paren{\tau_1\pP{p-k_1}(z)}-\dD F\Paren{0}}\vd \tau_1
    \\
    +
    \sum_{k_1=1}^p\alpha_{k_1}(s) z^{k_1}
    \int_0^1 \Paren[\Big]{\dD F\Paren{\tau_1 \pP{p}(z)}-
\dD F\Paren{\tau_1 \pP{p-k_1}(z)}}\vd \tau_1.
\end{multline}

Note that $\pP{p-k_1}(z)\in J(z;\balpha)$ for any $k_1 \in \{1,\ldots,p\}$ with $J(z;\balpha)$
defined in Notation~\ref{not:J}.
The Lipschitz constant of $\dD F$ on $J(z;\balpha)$ is the the sup-norm of $\dD^2 F$ on $J(z;\balpha)$, that is $\abs{F}_2(z)$ in Notation~\ref{not:J}.
 Then
\begin{multline*}
    \sum_{k_1=1}^p \abs{\alpha_{k_1}} \cdot\abs {z}^{k_1}\cdot
    \abs{\dD F\Paren{\tau_1\pP{p}(z)}-\dD F\Paren{\tau_1\pP{p-k_1}(z)}}
    \\
    \leq 
    \abs{F}_{2}(z) \sum_{k_1=1}^p \sum_{k_2=p-k_1+1}^{p}\abs{\alpha_{k_1}\alpha_{k_2}}\cdot \abs{z}^{k_1+k_2}
    \leq
    \abs{F}_{2}(z) \sum_{\substack{1\leq k_1,k_2\leq p
	\\
    k_1+k_2\geq p+1}}
    \abs{\alpha_{k_1}\alpha_{k_2}}\cdot \abs{z}^{k_1+k_2}.
\end{multline*}
This justifies our definition of $R_1(z)$ in Notation~\ref{not:remainder}
as $\abs{R_1(z)}=\grandO(\abs{z}^{p+1})$ for $z$ close to $0$.
Besides, 
\begin{equation*}
   \sum_{\substack{1\leq k_1,k_2\leq p\\ k_1+k_2\geq p+1}}
    \abs{\alpha_{k_1}\alpha_{k_2}}
    \leq \Paren*{\sum_{k=1}^p \abs{\alpha_k}}^2
    \leq \pP{p}^\star(1)^2.
\end{equation*}
This proves~\eqref{eq:ineq:Rm} for $R_1$.

We could then iterate the above decomposition in~\eqref{eq:100} by replacing $F$ by $\dD F$
for the second term in the right-hand side of the above integral. We
obtain
\begin{multline}
    \label{eq:3}
\dD F\Paren{\tau_1 \pP{p-k_1}(z)}-\dD F(0)
    \\
    = \tau_1 \dD^2 F(0) \pP{p-k_1}(z)
    +
    \sum_{\mathclap{k_2\leq p-k_1}}\tau_1\alpha_{k_2} z^{k_2}
\int_0^1 \Paren[\Big]{\dD^2 F\Paren{{\tau_2\tau_1}\pP{p-k_1-k_2}(z)}-\dD^2 F\Paren{0}}\vd \tau_2
    \\
    +
    \sum_{\mathclap{k_2 {\leq} p-k_1}}\tau_1\alpha_{k_2} z^{k_2}
\int_0^1 \Paren{ \dD^2 F\Paren[\Big]{\tau_1\tau_2 \pP{p-k_1}(z)}
    -
\dD^2 F\Paren{\tau_2\tau_1 \pP{p-k_1-k_2}(z) }}\vd \tau_2.
\end{multline}
Using a change of variable $\tau'_2=\tau_1\tau_2$  allows one to transform 
$\tau_1\int_0^1 f(\tau_1\tau_2)\vd \tau_2$ to $\int_0^{\tau_1} f(\tau_2)\vd\tau_2$
for any suitable function $f$. 

Integrating \eqref{eq:3} between $0$ and $1$ and using $p$ times the development \eqref{eq:100} to the successive
derivatives of $F$ leads to~\eqref{eq:102}.

The remainder's terms $R_k$ are controlled in a similar way as for $R_1$.
\end{proof}
%%%% end of proof

%%%%%%%%%%%%%%%%%%%%%%%%%%%%%%%%%%%%%%%%%%%%%%%%%%%%%%%%%%%%%%%%%%%%%%
%%%%%%%%%%%%%%%%%%%%%%%%%%%%%%%%%%%%%%%%%%%%%%%%%%%%%%%%%%%%%%%%%%%%%%
%%%%%%%%%%%%%%%%%%%%%%%%%%%%%%%%%%%%%%%%%%%%%%%%%%%%%%%%%%%%%%%%%%%%%%
\subsection{Application of the Taylor formula}

\label{sec:asy-inv-for}

In this section, we use Notations~\ref{not:delta},~\ref{not:remainder}-\ref{not:J}.

\begin{proposition}
    \label{prop:expansion}
	Assume Hypothesis~\ref{hyp:beta}.
	Let $F\in\cC^{p+1}(\Theta,\RR)^\bS$ and $s\in \bS$ 
	such that $\coef{F}_1(s) \neq 0$.
	Let $\balpha(s)\in \RR^{p+1}$ such that 
		$\eqref{eq:rs:s}$ is satisfied. 
    Assume that $\pP{m}(\varphi(s); \balpha(s)) \in \Theta$ for all $m \in \{1,\ldots,p\}$.
Then 
\begin{multline*}
     \abs{F(s,\pP{p}(\varphi(s);\balpha(s)))}
	 \leq 
	 \sum_{m \in \Set{2,\ldots, p+1}}
    \pP{p}^\star(1;\balpha(s))^{m} 
    \abs{F}_{m}(\varphi(s);\balpha(s))
    \abs{\varphi(s)}^{p+1}
	\\
	+ \sum_{\substack{m \in \Set{1,\ldots, p} \\ \beta_m+1 \leq \bbeta_\star} }
    \pP{p}^\star(1;\balpha(s))^{m}
    \abs{\coef{F}_m(s)}
    \cdot
    \abs{\varphi(s)}^{\gamma_m-\beta_m}.
\end{multline*}
\end{proposition}
\begin{proof}
In the proof we replace any dependence on $(\varphi(s);\balpha(s))$ by simply $(s')$.
Using Notation~\ref{not:delta}, we rewrite \eqref{eq:102} in Lemma~\ref{lem:devel:1}
as 
\begin{multline*}
    F(s,\pP{p}(s'))-\coef{F}_0(s) \varphi(s)^{-\beta_0}
    \\
    =
    \sum_{m=1}^p \frac{\coef{F}_m(s)}{m!}
    \sum_{k_1+\dotsb+k_m\leq p} \alpha_{k_1}(s)\dotsb\alpha_{k_{m}}(s)
    \varphi(s)^{k_1+\dotsb +k_{m}-\beta_m}+\sum_{m=1}^p R_m(s').
\end{multline*}
For $m=1,\dotsc,p$, the control \eqref{eq:ineq:Rm} implies that 
\begin{equation*}
    \abs{R_m(s')}\leq \abs{F}_{m+1}(s')\pP{p}^\star(1;\balpha(s))^{m+1}\abs{\varphi(s)}^{p+1}
    \text{ since }\abs{\varphi(s)}\leq 1
\end{equation*}
for all $m\in \{1,\ldots,p\}$.
Using the definition of $\gamma_m$ in~\eqref{eq:gamma}, we split 
\begin{multline*}
    \Set{(k_1,\dotsc,k_m)\given 1\leq k_1+\dotsb+k_m\leq p}
    \\
    =
    \Set{(k_1,\dotsc,k_m)\given 1 \leq k_1+\dotsb+k_m < \gamma_m}
    \\
    \cup
    \Set{(k_1,\dotsc,k_m)\given \gamma_m \leq k_1+\dotsb+k_m\leq p}.
\end{multline*}
Therefore, the convention that the sum is null if the index set is empty, 
yields for $m\in \{1,\ldots, p\}$ that if $\beta_m +1 \leq \bbeta_\star$ then
\begin{multline*}
    \abs*{\sum_{\gamma_m\leq k_1+\dotsb+k_m\leq p} \alpha_{k_1}(s)\dotsb\alpha_{k_{m}}(s)
    \abs{\varphi(s)}^{k_1+\dotsb +k_{m}}}
    \\
    \leq 
    \sum_{\gamma_m\leq k_1+\dotsb+k_m\leq p} 
    \abs{\alpha_{k_1}(s)\dotsb\alpha_{k_{m}}(s)}
    \cdot
    \abs{\varphi(s)}^{\gamma_m}
    \leq \pP{p}^\star(1;\balpha(s))^m\abs{\varphi(s)}^{\gamma_m}.
\end{multline*}
Using~\eqref{eq:rs:s} we obtain the result.
\end{proof}
%%
%%%%%%%%%%%%%%%%%%%%%%%%%%%%%%%%%%%%%%%%%%%%%%%%%%%%%%%%%%%%%%%%%%%%%%
\subsection{Proof of Theorem~\ref{thm:1} on the approximation of roots}

\begin{proof}[Proof of Theorem~\ref{thm:1}]
Since $\dD F(s,\cdot)$ is continuous on $U(s)$, 
The fact that $\abs{\dD F(s,\cdot)}$ has a positive lower bound on $U(s)$
(Condition~\ref{thm:root:iii})
implies that $F(s,\cdot)$ is invertible on $U(s)$
by the Inverse Function Theorem \cite{abraham}.
Besides, for any $\sigma,\sigma'\in F(s,U(s))$,
\begin{multline*}
    \abs{
    \dD F^{-1}(s,\sigma)-\dD F^{-1}(s,\sigma')
}
\leq \sup_{r\in[\sigma,\sigma']}\frac{\abs{\sigma-\sigma'}}{\abs{\dD F(s,F^{-1}(s,\tau))}}
\\
\leq \sup_{\theta\in U(s)} \frac{\abs{\sigma-\sigma'}}{\abs{\dD F(s,\theta)}}
= \sup_{\theta\in U(s)} \frac{\abs{\varphi(s)}^{\beta_1} }{\abs{\varphi(s)}^{\beta_1} \abs{\dD F(s,\theta)}} \abs{\sigma-\sigma'}
\leq \frac{\abs{\varphi(s)}^{\beta_1}}{c(s)}\abs{\sigma-\sigma'}.
\end{multline*}

Condition~\ref{thm:root:iv} implies that $J(\varphi(s))\subset U(s)$, 
where $J$ is defined in Notation~\ref{not:J}. 

With $g(s):=F(s,\pP{p}(\varphi(s);\balpha(s)))$ and since $F(s,\theta(s))=0$,  
\begin{equation*}
    \abs{\pP{p}(\varphi(s);\balpha(s))-\theta(s)}=\abs{F^{-1}(s,g(s))-F^{-1}(s,0)}
    \leq \frac{\abs{\varphi(s)}^{\beta_1}}{c(s)}\abs{g(s)}.
\end{equation*}
The result follows from Proposition~\ref{prop:expansion}.
\end{proof}
%%

%%%%%%%%%%%%%%%%%%%%%%%%%%%%%%%%%%%%%%%%%%%%%%%%%%%%%%%%%%%%%%%%%%%%%%

\printbibliography

\end{document}